\documentclass[a4,11pt]{article}

\usepackage[english]{babel}
\usepackage[utf8x]{inputenc}
\usepackage{amsmath,amssymb,amstext,amsthm}
\usepackage{graphicx}
\usepackage[colorinlistoftodos]{todonotes}
\usepackage[margin=1in]{geometry}
\usepackage{hyperref}
\usepackage{enumitem}

\newtheorem{theorem}{Theorem}[section]
\newtheorem{lemma}[theorem]{Lemma}
\newtheorem{corollary}[theorem]{Corollary}
\newtheorem{qn}[theorem]{Question}
\newtheorem{conj}[theorem]{Conjecture}

\theoremstyle{definition}
\newtheorem*{definition}{Definition}

\theoremstyle{definition}

\makeatletter
\renewcommand*{\@fnsymbol}[1]{\@arabic{#1}}
\makeatother

\newcommand{\ZZ}[1]{\mathbb{Z}_2[#1]}

\newcommand{\cart}{\mathbin{\square}}

\DeclareMathOperator{\aut}{Aut}

\DeclareMathOperator{\dist}{dist}

\newcommand{\qedn}{\qed\vspace{.2cm}}
\newcommand{\Ztn}{{\mathbb{Z}^n_2}}

\renewcommand{\geq}{\geqslant} 
 
\renewcommand{\ge}{\geqslant} 
\renewcommand{\le}{\leqslant}

\renewcommand{\restriction}{\mathord{\upharpoonright}}

\title{Cores of Cubelike Graphs}
\author{Laura Man\v{c}inska\thanks{QMATH, Department of Mathematical Sciences, University of Copenhagen, Denmark} \\
\texttt{mancinska@math.ku.dk}
\and
Irene Pivotto\thanks{School of Mathematics and Statistics, University of Western Australia, Australia} \\
\texttt{piv8irene@gmail.com}
\and
David E.~Roberson\thanks{Department of Physics, Technical University of Denmark, Denmark}\\
\texttt{dero@dtu.dk}
\and
Gordon F.~Royle\thanks{School of Mathematics and Statistics, University of Western Australia, Australia}\\
\texttt{gordon.royle@uwa.edu.au}}

\date{}

\begin{document}
\maketitle

\begin{abstract}
A graph is \emph{cubelike} if it is a Cayley graph for some elementary abelian $2$-group $\mathbb{Z}_2^n$. The core of a graph is its smallest subgraph to which it admits a homomorphism. More than ten years ago, Ne\v{s}et\v{r}il and \v{S}\'amal (On tension-continuous mappings. \emph{European J. Combin.,} 29(4):1025--1054, 2008) asked whether the core of a cubelike graph is cubelike, but since then very little progress has been made towards resolving the question. Here we investigate the structure of the core of a cubelike graph, deducing a variety of structural, spectral and group-theoretical properties that the core ``inherits'' from the host cubelike graph.  These properties constrain the structure of the core quite severely --- even if the core of a cubelike graph is not actually cubelike, it must bear a very close  resemblance to a cubelike graph. Moreover we prove the much stronger result that not only are these properties inherited by the core of a cubelike graph, but also by the orbital graphs of the core.  Even though the core and its orbital graphs look very much like cubelike graphs, we are unable to show that  this is sufficient to characterise cubelike graphs. However, our results are strong enough to eliminate all non-cubelike vertex-transitive graphs on up to $32$ vertices as potential cores of cubelike graphs (of any size). Thus, if one exists at all, a cubelike graph with a non-cubelike core has at least $128$ vertices and its core has at least $64$ vertices.
\end{abstract}

%\section{To Do}

%\input{intro}
\section{Introduction}

A \emph{homomorphism} from a graph $X$ to a graph $Y$ is an adjacency-preserving, but not necessarily injective, function from $V(X)$ to $V(Y)$. If such a function exists, then we say that $X$ \emph{maps} to $Y$, and write $X \rightarrow Y$.  As the pre-image of any vertex of $Y$ is a coclique of $X$ it is easy to see that a homomorphism $X \rightarrow K_k$ is equivalent to a $k$-coloring of $X$. Therefore the chromatic number of $X$ can be defined as the smallest $k$ such that $X \rightarrow K_k$. As a result, homomorphisms are often viewed as a generalization of colorings, and indeed a number of variants of graph coloring can be succinctly expressed in terms of the existence of homomorphisms to particular families of graphs. 

%For example, the \emph{fractional} chromatic number of a graph $X$ is the smallest rational number $v/k$ such that $X \rightarrow \kne{v}{k}$, where $\kne{v}{k}$ is the Kneser graph with $\binom{v}{k}$ vertices. % with vertex set the $k$-subsets of $\{1,2,\ldots,v\}$ and where two vertices are adjacent if and only if they are disjoint.

Two graphs $X$ and $Y$ are \emph{homomorphically equivalent}, denoted $X \leftrightarrow Y$, if  $X$ maps to $Y$ and $Y$ maps to $X$. Justifying its name, homomorphic equivalence is an equivalence relation on the set of all (unlabelled) graphs. As far as any  graphical parameter or property related to coloring is concerned, the graphs in each equivalence class are ``the same''. A graph $X$ is called a \emph{core} if it has no homomorphisms to any subgraph with fewer vertices (the name suggests some sort of ``incompressibility''). Each homomorphic equivalence class contains a unique core, which is thus a natural canonical representative for the equivalence class. If $X$ is a graph, then $X^\bullet$ denotes the core in the equivalence class containing $X$, and will be called \emph{the core} of $X$. 

%It is straightforward to see that $X^\bullet$ is an induced subgraph of $X$, and indeed for the graphs that we are considering, $X$ will have many induced subgraphs isomorphic to $X^\bullet$. Often we need to fix on a specific induced subgraph isomorphic to $X^\bullet$ and we indicate this by the phrase ``a copy of $X^\bullet$''. If $X$ is a graph, then we denote by 
%For example, any bipartite graph with at least one edge is homomorphically equivalent to $K_2$. A homomorphism from $X$ to itself is called an \emph{endomorphism} --- an endomorphism that is surjective is an automorphism, but otherwise it maps $X$ to a proper subgraph of itself. If $X$ has no endomorphisms other than automorphisms, then it is called a \emph{core} which can be loosely viewed as meaning that it is ``incompressible'' as far as coloring information is concerned.  Each homomorphism equivalence class contains a unique core, which is the graph with the fewest vertices in that class, and so the cores form a natural set of representatives for homomorphism equivalence classes. 
%A key tool used repeatedly in what follows is that if $Y$ is a copy of $X^\bullet$ then there is homomorphism $\varphi: X \rightarrow X$ with range equal to $Y$ and acting as the identity on $Y$. In other words, this homomorphism, which is called a \emph{retraction}, fixes $V(Y)$ vertex-wise and maps $X$ onto $Y$. We may say that a core $Y$ is a \emph{retract} of $X$ if it is the image of a retraction.

From the definitions, it follows that in general $X$ and its core will have the same values for any graph parameter that can be viewed as ``coloring-related'', such as the chromatic number or its variants like the fractional chromatic number.
% --- for example $\chi(X^\bullet) = \chi(X)$ and $\chi_{f}(X^\bullet)= \chi_f(X)$, where $\chi(X)$ and $\chi_f(X)$ are the \emph{chromatic number} and \emph{fractional chromatic number} of $X$ (respectively).  
More surprisingly however, is that various \emph{symmetry-related} properties of $X$ are shared by $X^\bullet$. Hahn and Tardif \cite{tardif} showed that if $X$ is vertex-transitive, then $X^\bullet$ is also vertex-transitive, and similarly for a number of stronger symmetry properties. 
%In addition, if the automorphism group of $X$ is primitive, then so is the automorphism group of $X^\bullet$. 
These connections are surprising because although the composition of an automorphism of $X$ and a retraction onto $X^\bullet$ is a mapping from $\mathrm{Aut}(X)$ to $\aut(X^\bullet)$, it is not a group homomorphism.  Although having a vertex-\emph{transitive} automorphism group is a symmetry property inherited by the core of a graph, the stronger property of having a vertex-\emph{regular} automorphism group is not, for there are Cayley graphs whose core is not a Cayley graph. In fact, 
any vertex-transitive core is the core of some Cayley graph, and so any non-Cayley vertex-transitive core provides an immediate example: the Petersen graph is an obvious choice here.

Although the class of all Cayley graphs is not closed under taking cores, there are particular families of Cayley graphs that might be. In this paper we consider \emph{cubelike  graphs}, which are the Cayley graphs for an elementary abelian group $\mathbb{Z}_2^n$, originally studied and named by Lov\'asz \cite{MR0398886} as a particular family of graphs with only integer eigenvalues.
%he most obvious example of a cubelike graph, and presumably the inspiration for the term ``cubelike'', is the $n$-dimensional cube $Q_n = \mathrm{Cay}(\mathbb{Z}_2^n, \{e_1, e_2, \ldots, e_n\})$, where $e_i$ is the $i$'th standard basis vector.  
Ne\v{s}et\v{r}il and \v{S}\'amal \cite{homotens} asked whether or not the core of a cubelike graph is cubelike, but given the supporting evidence, and the length of time that this question has been open, we feel that their question should be upgraded to a conjecture.

\begin{conj}
\label{conj:cubecore}
The core of a cubelike graph is cubelike.
\end{conj}

We consider the structure of the core of a cubelike graph, showing that such a graph is heavily constrained with respect to a variety of graphical, spectral and group-theoretical properties.  For graph-theoretical properties we show that if $Z$ is a cubelike graph with core $X$, then $X$ inherits the clique number and chromatic numbers of $Z$. If $Z$ is a cubelike graph of valency $d$, then there is a \emph{covering map} (a locally-injective surjective homomorphism) from the $d$-cube $Q_d \rightarrow Z$, which then implies that the eigenvalues of $Z$ are a sub-multiset of the eigenvalues of $Q_d$. We show that the core of a cubelike graph of valency $d$ is also covered by $Q_d$ and so its eigenvalues are also constrained.  For group-theoretical properties, we show the core of a cubelike graph has a generously-transitive automorphism group, which again is a property shared by all cubelike graphs.

Much more strongly we show that if $Z$ is a cubelike graph with core $X$, then it is not only $X$ that must have all of these properties, but all of the \emph{orbital graphs} of $X$ (that is, graphs with vertex set $V(X)$ whose automorphism group contains $\aut(X)$) must also share them. Given a particular core $X$, it can often be ruled out as a possible core for a cubelike graph by showing that it, or any one of its orbital graphs, does not have one or more of these properties.

However, despite these quite severe restrictions, we have not been able to show that they characterise cubelike graphs in general. To give some indication of how restrictive these conditions are, we apply the battery of tests to the $677116$ connected vertex-transitive graphs on $32$ vertices and demonstrate that the only graphs that meet all of our necessary conditions to be the core of a cubelike graph are themselves cubelike.  As a result, if a cubelike graph with a non-cubelike core exists, then it has at least $128$ vertices and its  core has at least $64$ vertices. We emphasize here that our results are looking ``upwards'' -- showing that a particular (non-cubelike) vertex-transitive graph cannot be the core of \emph{any} cubelike graph of any size, rather than the much easier situation of looking ``downwards'' and showing that a particular cubelike graph does in fact have a cubelike core. At least in principle, this approach could resolve the conjecture completely if enough additional properties of the core of a cubelike graph could be found that only a cubelike graph could satisfy them. 

In the study of homomorphisms of vertex-transitive graphs, a particular family of Cayley graphs seems to play a special role. We say that the Cayley graph $\text{Cay}(\Gamma, S)$ is \emph{normal} if $S$ is closed under conjugation by elements of $\Gamma$. Normal Cayley graphs have a number of attractive properties of particular relevance to coloring and other graph homomorphisms. In particular, Larose, Laviolette and Tardif \cite{LLT} showed that a graph $X$ is \emph{hom-idempotent}, i.e, homomorphically equivalent to its Cartesian square $X \cart X$, if and only if $X$ is homomorphically equivalent to a normal Cayley graph, and raised the following question.

\begin{qn}
\label{qn:normcore}
Is the core of a normal Cayley graph itself a normal Cayley graph?
\end{qn}

Every Cayley graph of an abelian group is a normal Cayley graph, but there are many vertex-transitive graphs that are not normal Cayley graphs.

\section{Preliminaries}\label{sec:preliminaries}

For us, a \emph{graph} $X$ consists of a vertex set $V(X)$ and edge set $E(X)$ of unordered pairs of distinct vertices of $X$. We indicate that two vertices $u$, $v$ are adjacent by $u \sim v$ and will use either $uv$ or $(u,v)$ to refer to the edge between $u$ and $v$ (which more formally would be $\{u,v\}$). While our main focus will be graphs, \emph{directed graphs} occur naturally in some parts of our work. A \emph{digraph} also consists of a vertex set and edge set, but each edge is now an ordered pair of vertices. We will sometimes use the term \emph{arc} to refer to a directed edge when the distinction is important. However, we will still use $u \sim v$ to indicate that there is an arc from $u$ to $v$, and similarly, $uv$ and $(u,v)$ will represent an arc directed from $u$ to $v$. For most of what follows, the context will determine when directions are important.
For graphs, the \emph{open neighborhood} of a vertex $u$, denoted $N(u)$, is the set of vertices adjacent to $u$. If $S$ is a set of (ordered) pairs of vertices of a graph $X$, then we will use $X(S)$ to denote the (di)graph on $V(X)$ with edge/arc set equal to $S$.

\paragraph{Homomorphisms.} A homomorphism $\varphi$ from a (di)graph $X$ to a (di)graph $Y$ is a function from $V(X)$ to $V(Y)$ such that $\varphi(u) \sim \varphi(v)$ whenever $u \sim v$. We refer to this property as ``preserving adjacency". A homomorphism from a (di)graph $X$ to itself is known as an \emph{endomorphism} of $X$. The \emph{image} of a homomorphism $\varphi$ from $X$ to $Y$ is the subgraph of $Y$ consisting of the vertices in $\varphi(V(X))$ and arcs/edges of $Y$ of the form $(\varphi(u),\varphi(v))$ where $u \sim v$. The set $\varphi(V(X)) \subseteq V(Y)$ will be referred to as the \emph{vertex image} of $\varphi$. The \emph{fibres} of a homomorphism $\varphi$ are the sets $\varphi^{-1}(u)$ where $u$ is in the vertex image of $\varphi$. Note that the adjacency-preserving property of homomorphisms implies that their fibres are independent sets (as $Y$ has no loops). A homomorphism $\varphi$ from $X$ to $Y$ is \emph{faithful} if its image is an \emph{induced} subgraph of $Y$. If a homomorphism $\varphi$ from $X$ to $Y$ is injective then its image is isomorphic to $X$, and thus the existence of an injective faithful homomorphism from $X$ to $Y$ is equivalent to $X$ being isomorphic to an induced subgraph of $Y$. %For any vertex $u$ in the vertex image of a homomorphism $\varphi$, the preimage $\varphi^{-1}(u)$ is called a \emph{fibre} of $\varphi$.

An isomorphism from a (di)graph $X$ to a (di)graph $Y$ is a bijection from $V(X)$ to $V(Y)$ such that $\varphi(u) \sim \varphi(v)$ if and only if $u \sim v$. We write $X \cong Y$ whenever there is an isomorphism from $X$ to $Y$. An isomorphism from $X$ to itself is an \emph{automorphism} of $X$. Note that if $\varphi: V(X) \to V(X)$ is a bijection, then to prove it is an automorphism of $X$ it suffices to show that it is an endomorphism, i.e., that it preserves adjacency. The automorphisms of a (di)graph $X$ form a permutation group known as the \emph{automorphism group of $X$}, which is denoted by $\aut(X)$.

Whenever there exists a homomorphism from $X$ to $Y$, we write $X \to Y$. It is not difficult to see that if $\varphi_1$ is a homomorphism from $X$ to $Y$ and $\varphi_2$ is a homomorphism from $Y$ to $Z$, then $\varphi_2 \circ \varphi_1$ is a homomorphism from $X$ to $Z$. Thus `$\to$' is a transitive relation, and it is trivially reflexive. It is not symmetric however,  since it is possible that $X \to Y$ but $Y \not\to X$, for instance if $X = K_n$ and $Y = K_{n+1}$. It is also not antisymmetric, as it is quite possible that $X \to Y$ and $Y \to X$, even when $X \ne Y$. If two graphs $X$ and $Y$ do satisfy  ``$X \to Y$ and $Y \to X$", then we write $X \leftrightarrow Y$ and say that $X$ and $Y$ are \emph{homomorphically equivalent}. It follows from the reflexivity and transitivity of `$\to$' that `$\leftrightarrow$' is indeed an equivalence relation on graphs. The relation `$\to$' induces a partial order on the equivalence classes of `$\leftrightarrow$', and this is known as the \emph{homomorphism order}.

\paragraph{Cores.} Each equivalence class of `$\leftrightarrow$' contains a unique graph with the fewest vertices and so these graphs form an obvious collection of ``canonical representatives'' for the equivalence relation.  These graphs are known as \emph{cores} and are the primary focus of our work, so we derive some of their fundamental properties here. Most of these facts are straightforward to prove, but we refer the reader to~\cite{AGT} or~\cite{tardif} for full explanations.

A (di)graph $X$ is a \emph{core} if all of its endomorphisms are bijections (and therefore automorphisms). A non-bijective endomorphism is referred to as a \emph{proper} endomorphism, so a core can also be defined as a graph with no proper endomorphisms. A (di)graph $X$ is a \emph{core of a (di)graph $Y$} if $X$ and $Y$ are homomorphically equivalent and $X$ is a core.
Every (di)graph $X$ has a unique core up to isomorphism, and we denote this by $X^\bullet$ and refer to it as \emph{the} core of $X$. The core of $X$ is always an induced subgraph of $X$, and in fact can be defined as the vertex-minimal endomorphic image of $X$. Two (di)graphs $X$ and $Y$ are homomorphically equivalent if and only if $X^\bullet \cong Y^\bullet$. The core of $X$ is the unique graph with the fewest vertices in the homomorphic equivalence class of $X$.

Finally, if $X'$ is a subgraph of $X$ isomorphic to its core, then there exists an endomorphism $\varphi$ of $X$ with image $X'$ such that $\varphi(x) = x$ for all $x \in V(X')$, i.e., $\varphi^2 = \varphi$. Such an endomorphism, which acts as the identity on its image, is called a \emph{retraction} and its image a \emph{retract}.

\paragraph{Colorings, cliques, odd cycles, etc.}

A \emph{$c$-coloring} of a graph $X$ is an assignment of ``colors" from the set $[c] := \{1, \ldots, c\}$ to the vertices of $X$ such that adjacent vertices receive different colors. It straightforward to see that a $c$-coloring is equivalent to a homomorphism to the complete graph $K_c$. Thus the \emph{chromatic number} of $X$, defined as the smallest number of colors required to color $X$, is equal to the size of the smallest complete graph to which $X$ admits a homomorphism. Many other graph parameters, especially ``coloring-related'' parameters, can be defined in terms of homomorphisms. Two in particular that we will use are the clique number and odd girth of a graph $X$, denoted $\omega(X)$ and $\text{og}(X)$ respectively. The \emph{clique number} is the size of a largest complete subgraph of $X$, and this is equal to $\max\{r : K_r \to X\}$. The odd girth of $X$ is the length of a shortest odd cycle of $X$ (defined to be $\infty$ if $X$ is bipartite). It is not difficult to show that $\text{og}(X) = \min\{g: C_g \to X, g \text{ odd}\}$, where $C_g$ is the cycle of length $g$. It follows immediately from these reformulations in terms of homomorphisms that if $X \to Y$, then $\omega(X) \le \omega(Y)$, $\chi(X) \le \chi(Y)$, and $\text{og}(X) \ge \text{og}(Y)$. Thus, if $X$ and $Y$ are homomorphically equivalent, then these three parameters, and any others that can be similarly reformulated, must be equal for $X$ and $Y$. In particular, this implies that such parameters are the same for $X$ and $X^\bullet$.

Note that $\omega(X) \le \chi(X)$ for any graph $X$ and, moreover, these two parameters are equal (to $r$) if and only if $X$ is homomorphically equivalent to a complete graph (of size $r$). This implies that the core of $X$ is a complete graph if and only if its clique and chromatic numbers are equal.

Another parameter that will be important to us is the \emph{independence number} of $X$, which is the maximum number of pairwise non-adjacent vertices in $X$, and is denoted $\alpha(X)$. Note that $\alpha(X) = \omega(\overline{X})$, where $\overline{X}$ is the \emph{complement} of $X$. Unlike the parameters mentioned above, the independence number is not monotone with respect to homomorphisms, i.e., $X \to Y$ does not immediately tell us anything about how $\alpha(X)$ and $\alpha(Y)$ compare. Indeed, adding many isolated vertices to either $X$ or $Y$ does not change whether $X \to Y$, but it allows us to increase the independence number of either arbitrarily. On the other hand, for a vertex-transitive graph $X$, it is known that the fractional chromatic number is equal to $|V(X)|/\alpha(X)$~\cite{AGT}. Since the fractional chromatic number is monotone with respect to homomorphisms, this means that if $X$ and $Y$ are vertex-transitive graphs, then $X \to Y$ implies that $|V(X)|/\alpha(X) \le |V(Y)|/\alpha(Y)$.

\paragraph{Covering maps and eigenvalues.}

Let $\varphi$ be a homomorphism from a graph $X$ to a graph $Y$. Then $\varphi$ is a \emph{covering map} if it is surjective and \emph{locally bijective}. The latter means that for every $u \in V(X)$, the restriction of the map $\varphi$ to the neighborhood of $u$ is a bijection between $N(u)$ and $N(\varphi(u))$. This implies that $u$ and $\varphi(u)$ have the same degree, and that $\varphi$ also acts surjectively on edges. If such a map exists from $X$ to $Y$, we say that $X$ \emph{covers} $Y$ or that $X$ is a \emph{cover} of $Y$.

Given a graph $X$, its \emph{adjacency matrix} $A$ is the symmetric $01$-matrix with rows/columns indexed by $V(X)$ such that $A_{uv} = 1$ if $u \sim v$ and $A_{uv} = 0$ otherwise. The eigenvalues of a graph are then the eigenvalues of its adjacency matrix. Since it is symmetric, it has a full set of orthonormal eigenvectors that span $\mathbb{R}^{V(X)}$. It is often useful to think of an eigenvector of $X$ as a function on its vertices, and we shall do so wherever convenient.

%A partition $\mathcal{P} = \{P_1, \ldots, P_r\}$ of the vertex set of a graph $X$ is said to be \emph{equitable} if there exist numbers $b_{ij}$ such that for any vertex $u \in P_i$, the number of neighbors of $u$ in $P_j$ is equal to $b_{ij}$. Equivalently, the parts $P_i$ induce regular subgraphs of $X$ and the edges between two distinct parts form a semiregular bipartite subgraph. The $r \times r$ matrix $B = (b_{ij})$ is said to be the \emph{quotient matrix} of $\mathcal{P}$. Brouwer and Haemers~\cite{spectraofgraphs}[Lemma~2.3.1], show that the eigenvalues of $B$ are a submultiset of the eigenvalues of $X$.

Suppose that $\varphi$ is a covering map from $X$ to $Y$. If $f :V(Y) \to \mathbb{C}$ is an eigenvector of $Y$, then $f \circ \varphi: V(X) \to \mathbb{C}$ is an eigenvector of $X$ with the same eigenvalue. (see Brouwer and Haemers \cite[Section 6.4]{spectraofgraphs}). It is easy to see that the linear map $f \mapsto f \circ \varphi$ is injective. Therefore, if $X$ covers $Y$, then the eigenvalues of $Y$ are a submultiset of the eigenvalues of $X$. This will be important for Corollary~\ref{cor:spectra}.

%Suppose that $\varphi$ is a covering map from $X$ to $Y$. We will show that the fibres of $\varphi$ form an equitable partition of $X$. First, the fibres are independent sets, i.e., they induce $0$-regular subgraphs of $X$. Now consider distinct $u,v \in V(Y)$. If $u \not\sim v$, then there are no edges between $\varphi^{-1}(u)$ and $\varphi^{-1}(v)$ since $\varphi$ preserves adjacency. Thus, in this case every vertex of $\varphi^{-1}(u)$ has exactly zero neighbors in $\varphi^{-1}(v)$. Now suppose that $u \sim v$. If $x \in \varphi^{-1}(u)$ then $\varphi$ maps the neighbors of $x$ bijectively to the neighbors of $u$. Therefore, $x$ has exactly one neighbor in $\varphi^{-1}(v)$. This proves that the partition of $V(X)$ into the fibres of $\varphi$ is equitable. Moreover, if $B$ is the quotient matrix of this partition, then the above argument shows that $B_{uv} = 1$ if $u \sim v$, and $B_{uv} = 0$ otherwise. In other words, $B$ is the adjacency matrix of $Y$. Putting this all together, we see that if $X$ covers $Y$, then the eigenvalues of $Y$ are a submultiset of the eigenvalues of $X$.

\subsection{Cayley Graphs}

Let $\Gamma$ be a group and let $C \subseteq \Gamma$. The \emph{Cayley digraph} on $\Gamma$ with \emph{connection set} $C$, denoted $\text{Cay}(\Gamma,C)$, has vertex set $\Gamma$ and there is an arc from $a$ to $b$ if $ba^{-1} \in C$. In other words, the arcs of $\text{Cay}(\Gamma,C)$ have the form $(a,ca)$ for $a \in \Gamma$ and $c \in C$. If $C$ does not contain the identity and is inverse closed (i.e., $c \in C \Leftrightarrow c^{-1} \in C$), then $\text{Cay}(\Gamma,C)$ has no loops and there is an arc from $a$ to $b$ if and only if there is an arc from $b$ to $a$. In this case we consider $\text{Cay}(\Gamma,C)$ as an undirected graph.

Given a Cayley (di)graph $X = \text{Cay}(\Gamma,C)$, note that right multiplication by an element of $\Gamma$ is an automorphism of $X$. Indeed, $a \sim b \Leftrightarrow ba^{-1} \in C \Leftrightarrow bg(ag)^{-1} \in C \Leftrightarrow ag \sim bg$. It quickly follows that $\Gamma$ is a subgroup of $\aut(X)$ that acts regularly on $V(X)$.
%Note that for any pair of elements $a,b \in \Gamma$, there exists a unique element of $\Gamma$, namely $a^{-1}b$, such that right multiplication by that element maps $a$ to $b$. Thus $\Gamma$ is a subgroup of $\aut(X)$ that acts regularly on $V(X)$.
Conversely, a (di)graph $X$ is a Cayley (di)graph if and only if $\aut(X)$ contains a regular subgroup. Note that the action of $\Gamma$ on $\text{Cay}(\Gamma,C)$ implies that it is \emph{vertex transitive}, i.e., for any two vertices $u,v$ there exists an automorphism mapping $u$ to $v$.

For the remainder of this section we will only consider the case where the connection set $C$ is inverse closed and thus $\text{Cay}(\Gamma,C)$ is a Cayley graph. We remark that if $\sigma$ is a group automorphism  of $\Gamma$, then $g \mapsto \sigma(g)$ is a (di)graph isomorphism from $\text{Cay}(\Gamma, C)$ to $\text{Cay}(\Gamma, \sigma(C))$.

\paragraph{Normal Cayley graphs.} Unlike right multiplication, left multiplication by group elements is not necessarily an automorphism, since it is possible that $ba^{-1} \in C$ but $(gb)(ga)^{-1} = g(ba^{-1})g^{-1} \not\in C$. In fact, it is easy to see that left multiplication by an element $g$ is an automorphism of the Cayley graph $X = \text{Cay}(\Gamma,C)$ for all $g \in \Gamma$ if and only if $C$ is closed under conjugation by any group element of $\Gamma$, i.e., $C = gCg^{-1} := \{gcg^{-1} : c \in C\}$. This is equivalent to $C$ being the union of conjugacy classes of $\Gamma$. In this case, $X = \text{Cay}(\Gamma,C)$ is called a \emph{normal Cayley graph}\footnote{The term \emph{normal Cayley graph} is also used to refer to Cayley graphs $X$ on a group $\Gamma$ where the subgroup of $\aut(X)$ corresponding to right multiplication by elements of $\Gamma$ is a normal subgroup of $\aut(X)$.}. Note that if $\Gamma$ is abelian, then $X = \text{Cay}(\Gamma,C)$ is always normal.

We saw above that Cayley graphs are vertex transitive. It turns out that normal Cayley graphs are additionally \emph{generously transitive}, i.e., for any two vertices $u$ and $v$ there is an automorphism that maps $u$ to $v$ and $v$ to $u$. This is well known, but we give a proof for completeness.

\begin{lemma}\label{lem:gentrans}
Every normal Cayley graph is generously transitive.
\end{lemma}
\proof
Let $\Gamma$ be a group and $C \subseteq \Gamma \setminus \{1\}$ be a union of conjugacy classes of $\Gamma$ so that $X = \text{Cay}(\Gamma,C)$ is a normal Cayley graph. We first show that the inverse map $\iota: \Gamma \to \Gamma$ given by $\iota(g) = g^{-1}$ is an automorphism of $X$. Obviously, $\iota$ is a bijection, so we merely need to show that it preserves adjacency. Suppose that $a \sim b$, i.e., $ba^{-1} \in C$. Since $C$ is closed under conjugation and inverses, we have that $b^{-1}a = b^{-1}(ab^{-1})b = b^{-1}(ba^{-1})^{-1}b\in C$. Therefore, $a^{-1} \sim b^{-1}$ and so $\iota$ is an automorphism of $X$.

Now let $a,b \in \Gamma$. Let $\sigma_{a^{-1}}$ and $\sigma_b$ be the automorphisms of $X$ given by right multiplication by $a^{-1}$ and $b$ respectively. Consider the automorphism $\sigma_b \circ \iota \circ \sigma_{a^{-1}}$. We have that
\[\sigma_b \circ \iota \circ \sigma_{a^{-1}} (a) = \sigma_b \circ \iota (1) = \sigma_b(1) = b,\]
and 
\[\sigma_b \circ \iota \circ \sigma_{a^{-1}} (b) = \sigma_b \circ \iota (ba^{-1}) = \sigma_b(ab^{-1}) = a.\]\qedn

%c = db <=> d^-1 = bc^-1 <=> c^-1d^-1c = c^-1b <=> (c^-1dc)^-1 = c^-1b
%
%suppose c = db and x = yb and c^-1dc = x^-1yx
%
%then b^-1d^-1ddb = b^-1y^-1yyb
%and so b^-1db = b^-1yb
%and so d = y
%and so c = db = yb = x
%
%
%
%c^-1 =  d^-1b <=> dc^-1 = b
%
%(gcg^-1)^-1a = gc^-1g^-1a
%
%
%
%000000000
%021000000
%
%210000000
%201000000
%
%120000000

%Suppose that $X = \text{Cay}(\Gamma,C)$ is a normal Cayley graph. Consider the inverse map on $\Gamma$, i.e., $g \mapsto g^{-1}$. Obviously, this is a bijection. Furthermore, if $a \sim b$ then $ab^{-1} \in C$ and thus $b^{-1}a = a^{-1}(ab^{-1})a \in C$ since $C$ is closed under conjugation. This implies that $a^{-1} \sim b^{-1}$, and so we see that the inverse map is always an automorphism of a normal Cayley graph (in fact, it is known that a Cayley graph is normal if and only if the inverse map is an automorphism). Using this, we can show that a normal Cayley graph $X$ is \emph{generously transitive}, i.e., for any two vertices $u,v$ there exists an automorphism $\sigma$ such that $\sigma(u) = v$ and $\sigma(v) = u$. To see this recall that $X$ is vertex transitive since it is a Cayley graph, and so it suffices to show that for any element $a \in \Gamma$, there is an automorphism of $X$ swapping $1$ and $a$. For this, we compose the inverse map with right multiplication by $a$. This maps 1 to $1^{-1}a = a$ and maps $a$ to $a^{-1}a = 1$ as desired.

\paragraph{The clique-coclique bound.}

A coclique is another name for an independent set, and the clique-coclique bound states that if $X$ is vertex transitive, then $\alpha(X) \omega(X) \le |V(X)|$. Note that for any graph $X$ we have $\chi(X) \ge |V(X)|/\alpha(X)$, since color classes must be independent sets. If $X$ is vertex transitive and $\omega(X) = \chi(X)$, then by the clique-coclique bound we have that $|V(X)|/\alpha(X) \ge \omega(X) = \chi(X) \ge |V(X)|/\alpha(X)$. Thus we have equality throughout and so the clique-coclique bound must hold with equality. The converse does not hold in general, even for Cayley graphs. For instance, let $X$ be the the graph obtained from the 1-skeleton of the cuboctahedron by adding edges between pairs of vertices at distance three. Then $X$ is a Cayley graph on 12 vertices with $\omega(X) = 3$ and $\alpha(X) = 4$, but $\chi(X) = 4$.

However, it is known~\cite{chrisnotes} that if $X$ is a normal Cayley graph, then the converse does hold. Thus for a normal Cayley graph $X$ the clique-coclique bound holds with equality if and only if $\omega(X) = \chi(X)$ if and only if $X$ has a complete graph as a core. Using the notion of homomorphic equivalence, we can extend this slightly to the following:

\begin{lemma}\label{lem:cliquecoclique}
Let $X$ be a vertex-transitive graph that is homomorphically equivalent to a normal Cayley graph. Then $\omega(X) = \chi(X)$ if and only if $\alpha(X)\omega(X) = |V(X)|$.
\end{lemma}
\proof
Let $Y$ be a normal Cayley graph homomorphically equivalent to $X$.  
Since $\omega(X) = \omega(Y)$ and $\chi(X) = \chi(Y)$, these two parameters are equal for $X$ if and only if they are equal for $Y$. Also, since they are vertex transitive and homomorphically equivalent, we have that $|V(X)|/\alpha(X) = \chi_f(X) = \chi_f(Y) = |V(Y)|/\alpha(Y)$. Therefore the clique-coclique bound holds with equality for $X$ if and only if it does so for $Y$. Since the statement of the lemma holds for normal Cayley graphs, the above two equivalences prove the statement for $X$.\qedn

In Section~\ref{subsec:prevresults}, we will see that the core of a vertex-transitive graph must be vertex transitive. Combining this with the above observations provides a useful tool for ruling out potential cores of normal Cayley graphs.

\paragraph{Cubelike graphs.}

The objects that are the main focus of this work are the Cayley graphs for $\mathbb{Z}_2^n$ for some $n \in \mathbb{N}$. These are known as \emph{cubelike graphs} since they are generalizations of the $n$-cube graph $Q_n$, which is $\text{Cay}(\mathbb{Z}_2^n,\{e_1, \ldots, e_n\})$, where $e_i$ is the \emph{$i^\text{th}$ standard basis vector} of $\mathbb{Z}_2^n$. This group is also a vector space and we will often take this view, considering the elements of $\mathbb{Z}_2^n$ as vectors, which in turn we may view as binary strings. For this group $-a = a$, and thus $a - b = a+b$, and so we will use the latter throughout. We will refer to the \emph{weight} of an element of $\mathbb{Z}_2^n$ as the number of 1's appearing in its representation as a vector/binary string. We remark that the automorphism group of $\mathbb{Z}_2^n$ is the group $\text{GL}(n,2)$ of invertible linear maps. If a cubelike graph $Z = \text{Cay}(\mathbb{Z}_2^n,C)$ is connected, then $C$ must contain a basis of $\mathbb{Z}_2^n$, and so we may assume without loss of generality that it contains all of the standard basis vectors. As $\mathbb{Z}_2^n$ is abelian, all of its subsets are closed under conjugation, and so cubelike graphs are always normal Cayley graphs. There are two special families of cubelike graphs that we will need later.

\begin{itemize}
\item \textbf{Folded cube:} For a given $n \in \mathbb{N}$ with $n > 1$, let $C = \{e_1, \ldots, e_{n-1}, e_1 + \ldots +e_{n-1}\}$. Then the cubelike graph $\text{Cay}(\mathbb{Z}_2^{n-1},C)$ is known as the \emph{folded cube of order $n$}. This graph can be constructed from the $(n-1)$-cube by adding edges between vertices at maximum distance $n-1$ from each other, or by identifying pairs of vertices of the $n$-cube at maximum distance from each other. The order $n$ of the folded cube coincides with its degree. For $n$ even this graph is bipartite, and for $n$ odd it has chromatic number 4. This graph is always distance transitive, i.e, for any two pairs of vertices $(a_1,b_1)$ and $(a_2,b_2)$ such that the distance between $a_i$ and $b_i$ is the same for $i = 1,2$, there is an automorphism mapping the first pair to the second. There is exactly one dependency among the vectors of $C$, namely that their sum is 0. This property of the connection set characterizes the folded cube. From this it is not too difficult to see that for $n$ odd the odd girth of the folded cube is $n$, and every pair of vertices is contained in a shortest odd cycle. For $n = 5$, the folded cube is better known as the \emph{Clebsch graph}. 

\item \textbf{Halved cube:} For a given $n \in \mathbb{N}$ with $n > 1$,  let $C = \{e_i : i \in [n-1]\} \cup \{e_i+e_j : i,j \in [n-1], i \ne j\}$. The cubelike graph $\text{Cay}(\mathbb{Z}_2^{n-1},C)$ is the \emph{halved $n$-cube}, denoted $\frac{1}{2}Q_{n}$. The distance-2 graph of the $n$-cube has two connected components  (the even / odd weight vertices respectively), and each is isomorphic to $\frac{1}{2}Q_{n}$. It is also isomorphic to the graph obtained from $Q_{n-1}$ by adding edges between vertices at distance two. The clique number of $\frac{1}{2}Q_n$ is $n$ unless $n = 3$ in which case the clique number is $4$. For $n=5$, the halved cube is the complement of the Clebsch graph.
\end{itemize}

%clique and chromatic number not equal to 3

\section{First Observations}

Here we discuss some of the basic properties of endomorphisms and cubelike graphs that we will make use of throughout the rest of the paper. We begin by presenting some previously known results about cores of cubelike graphs. Next we consider two ideas for proving Conjecture~\ref{conj:cubecore}, and then show why these ideas do not work. Finally, we prove some new results about endomorphisms of graphs that we will use for our later results.

\subsection{Previous results}\label{subsec:prevresults}

We discussed above that the core of a graph $X$ ``inherits" the homomorphic properties of $X$, such as clique and chromatic numbers. This means that if there is some value that these parameters cannot attain on cubelike graphs, then this value can also not be attained on cores of cubelike graphs. The first and easiest example is clique number. If a cubelike graph contains a clique of size three, then it must contain a clique of size four. Indeed, if the vertices $a,b,c$ form a triangle in a cubelike graph, then the vertices $a,b,c,a+b+c$ form a $K_4$. Therefore $\omega(X) \ne 3$ for any cubelike graph $X$, and thus the core of a cubelike graph cannot have clique number 3 either.

A similar, but nontrivial, result is known about the chromatic number of cubelike graphs. Payan~\cite{payan} proved that the chromatic number of a cubelike graph can never be equal to three. Thus the core of a cubelike graph cannot have chromatic number equal to three either. These two results will be useful when ruling out potential cores of cubelike graphs in Section~\ref{subsec:32vertices}.

It turns out that if $X$ has certain symmetry properties, then these are also inherited by $X^\bullet$. For instance, it is known that if $X$ is vertex transitive, then so is $X^\bullet$~\cite{AGT,chrisnotes,tardif}. This holds for many other types of symmetry as well, in particular it holds for generous transitivity. Since any normal Cayley graph is generously transitive by Lemma~\ref{lem:gentrans}, the core of a normal Cayley graph must be generously transitive. This applies to the core of cubelike graphs as well, since cubelike graphs are normal Cayley graphs.

Another property of the automorphism group of a graph that is inherited by cores is that of primitivity. A permutation group is said to be \emph{primitive} if it preserves no nontrivial partition of the set it acts on. It has been shown that if $X$ has a primitive automorphism group, then so does the core of $X$, though the proof of this is more difficult than the proof for vertex transitivity. This result appears in~\cite{chrisnotes} where it is credited to Tardif.

Finally, it is a result of Hahn and Tardif~\cite{tardif} that if $X$ is vertex transitive, then $|V(X^\bullet)|$ divides $|V(X)|$. Since cubelike graphs are vertex transitive and have $2^n$ vertices for some $n$, this result implies that the core of a cubelike graph has $2^k$ vertices for some $k$.

Summarizing all of the above, we see that the core of a cubelike graph:
\begin{itemize}[noitemsep,topsep=6pt]
\item has $2^k$ vertices for some $k \in \mathbb{N}$,
\item is generously transitive, and
\item does not have clique or chromatic number equal to three.
\end{itemize}

\subsection{Two false leads}\label{subsec:falseleads}

%\lnote{Probably not the best section title. Any suggestions?}
%\gnote{Is this better?}

At first sight, there does not seem to be any strong reason to believe that Conjecture~\ref{conj:cubecore} should hold. Although the core of a vertex-transitive graph is vertex transitive, there are many Cayley graphs whose cores are not Cayley graphs. In fact, it is known that \emph{any} vertex-transitive core is the core of \emph{some} Cayley graph, and there are many known non-Cayley cores, for example the Kneser graphs. Since a graph $X$ is a Cayley graph for a group $\Gamma$ if and only if $\Gamma$ is a regular subgroup of $\aut(X)$, this means that the property of $\aut(X)$ containing a regular subgroup is not preserved by taking cores. Thus Conjecture~\ref{conj:cubecore} asserts that something very special happens when the subgroup in question is $\mathbb{Z}_2^n$.

In our experience, researchers encountering this question for the first time almost always quickly arrive at one (or both) of two stronger statements, either of which would imply that cores of cubelike graphs are cubelike. Unfortunately, although each of these statements is plausible, neither is true. 
Each of the statements essentially says that there is a copy of the core embedded ``nicely'' in the host cubelike graph --- sufficiently nicely that it can rapidly be shown that it is indeed cubelike. We explain these two ideas and give counterexamples to each in order to illustrate some the difficulties involved in precisely identifying a copy of the core, and to justify the lengths to which we are forced to go to gather information about the core. If this helps dissuade future researchers from going down these same dead-ends, then that will be a bonus.

\paragraph{False lead 1: There is a copy of the core on a subgroup of vertices.} Consider a Cayley graph $X = \text{Cay}(\Gamma,C)$. If $Y$ is an induced subgraph of $X$ such that $V(Y) = \Gamma'$ where $\Gamma'$ is a subgroup of $\Gamma$, then it is easy to see that $Y$ is isomorphic to $\text{Cay}(\Gamma',C \cap \Gamma')$ and is therefore a Cayley graph for the group $\Gamma'$. Since any subgroup of $\mathbb{Z}_2^n$ is isomorphic to $\mathbb{Z}_2^k$ for some $k$, if a cubelike graph $Z$ contains a copy of its core on a subgroup of vertices, then its core must be cubelike. Unfortunately this does not always happen for cubelike graphs. Consider the halved 8-cube $\frac{1}{2}Q_8$. The set $\{0,e_1, \ldots, e_7\}$ induces a $K_8$ in $\frac{1}{2}Q_8$, and we will see that this is its core. Moreover, it is not hard to see that the largest subgroup of vertices forming a clique has size 4. Indeed, any such subgroup must contain only neighbors of $0$, i.e., the weight one and weight two elements. The largest subgroup containing only elements of weight at most one has size two. If a subgroup contains only elements of weight at most two, then every nonzero element must have a 1 in common with every weight two element. This leaves $\{0,e_1,e_2, e_1+e_2\}$ and $\{0,e_1+e_2,e_2+e_3,e_1+e_3\}$ as the only two possibilities up to permutation of the coordinates. Thus there is no copy of the core of $\frac{1}{2}Q_8$ on a subgroup of vertices.

To show that $K_8$ is the core of $\frac{1}{2}Q_8$, we must show that we can color the graph with 8 colors. To do this, fix a bijective assignment $f: [7] \to \mathbb{Z}_2^3 \setminus {0}$, of the nonzero elements of $\mathbb{Z}_2^3$ to the coordinates of the elements of $\mathbb{Z}_2^7 = V(\frac{1}{2}Q_8)$. We will then color the vertices of $\frac{1}{2}Q_8$ with the 8 elements of $\mathbb{Z}_2^3$ according to the function $\varphi: \mathbb{Z}_2^7 \to \mathbb{Z}_2^3$ given by
\[\varphi\left(\sum_{i \in S \subseteq [7]} e_i\right) = \sum_{i \in S} f(i).\]
Note that this map is linear, i.e., it is a group homomorphism. If $a \sim b$ in $\frac{1}{2}Q_8$, then $a+b = e_i$ or $a+b = e_i+e_j$ for some $i \ne j \in [7]$. Thus $\varphi(a) + \varphi(b) = f(i) \ne 0$ in the former case, or $\varphi(a)+\varphi(b) = f(i) + f(j) \ne 0$ in the latter case. In either case, $\varphi(a) \ne \varphi(b)$, and so this is a proper coloring of $\frac{1}{2}Q_8$.

So we have shown that the core of $\frac{1}{2}Q_8$ is $K_8$, but that $\frac{1}{2}Q_8$ contains no copy of $K_8$ on a subgroup of vertices. Note that the same argument can be used to show that $\frac{1}{2}Q_{2^r}$ has core $K_{2^r}$ but contains no copy of this on a subgroup of vertices, as long as $r \ge 3$.

\paragraph{False lead 2: There is a homomorphism to the core whose fibres are the cosets of some subgroup.}
In Lemma~\ref{lem:cosets}, we show that if $Z$ is a Cayley graph on $\mathbb{Z}_2^n$ and $\varphi$ is a vertex- and edge-surjective homomorphism to a graph $X$ such that the fibres of $\varphi$ are cosets of some fixed subgroup of $\mathbb{Z}_2^n$, then $X$ must be a cubelike graph. Any homomorphism from a graph to its core must be both vertex- and edge-surjective. Thus if a cubelike graph has a homomorphism to its core whose fibres are cosets of some fixed subgroup, then its core must also be cubelike. Unfortunately, it is not always possible to find such a homomorphism to the core of a cubelike graph. The counterexample here is actually the complement of the counterexample above.

Let $Z = \frac{1}{2}Q_8$ be the halved 8-cube. Recall from above that we showed that the core of $Z$ is $K_8$. Thus $\chi(Z) = \omega(Z) = 8$. By Lemma~\ref{lem:cliquecoclique}, this implies that $\alpha(Z)\omega(Z) = |V(Z)|$. Therefore, the complement, $\overline{Z}$, satisfies $\omega(\overline{Z})\alpha(\overline{Z}) = |V(\overline{Z})|$. Applying Lemma~\ref{lem:cliquecoclique} again, we have that
\[\chi(\overline{Z}) = \omega(\overline{Z}) = |V(\overline{Z})|/\alpha(\overline{Z}) =  |V(Z)|/\omega(Z) = 128/8 = 16.\]
Note that any homomorphism from $\overline{Z}$ to its core is a 16-coloring of $\overline{Z}$, and so the possible fibres of a homomorphism to its core are the same as the possible color classes of a 16-coloring. Since $\chi(\overline{Z}) = |V(\overline{Z})|/\alpha(\overline{Z})$, the color classes in a 16-coloring of $\overline{Z}$ must be maximum independent sets of $\overline{Z}$, i.e., maximum cliques of $Z$. However, above we showed that there are no maximum cliques of $Z$ consisting of vertices forming a subgroup. Therefore, the color classes of a 16-coloring of $\overline{Z}$ cannot be cosets of a fixed subgroup, and thus neither can the fibres of a homomorphism from $\overline{Z}$ to its core. In fact, not a single fibre can be a coset of some subgroup.

Note that the same argument can be used to show that, for $ r \ge 3$, the core of the complement of $\frac{1}{2}Q_{2^r}$ is $K_{2^{2^r - r - 1}}$, but no homomorphism to the core has fibres that are cosets of a fixed subgroup.\\

These examples show that we cannot prove Conjecture~\ref{conj:cubecore} via either of these statements. This is somewhat discouraging since each of these plausible (but unfortunately false) statements would have provided a satisfactory explanation for \emph{why} Conjecture~\ref{conj:cubecore} should be true.

\subsection{Properties of Endomorphisms}

In this section, we prove some preliminary lemmas that we will need later, and that may be useful more generally. The basic idea is that endomorphisms of a graph $Z$ must act as isomorphisms between the subgraphs of $Z$ isomorphic to its core. More precisely, we have the following:

\begin{lemma}\label{lem:iso}
Let $Z$ be a graph with an induced subgraph $X$ which is isomorphic to the core of $Z$. If $\varphi$ is an endomorphism of $Z$ then, when restricted to $V(X)$, the map $\varphi$ acts as an isomorphism from $X$ to the graph induced by $\varphi(V(X))$.
\end{lemma}
\proof
Let $Y$ be the subgraph induced by $\varphi(V(X))$ and let $\rho$ be a retraction of $Z$ onto $X$. Obviously, the restriction $\varphi\restriction_X$ is a homomorphism from $X$ to $Y$ and the restriction $\rho\restriction_Y$ is a homomorphism from $Y$ to $X$. Thus $\rho\restriction_Y \circ \varphi\restriction_X$ is a homomorphism from $X$ to itself and is therefore an automorphism since $X$ is a core. It follows that $\varphi_X$ preserves non-edges (thus is injective), and it is surjective by definition. Therefore $\varphi\restriction_X$ is an isomorphism from $X$ to $Y$.\qedn

%Obviously, $\varphi$ acts injectively on $V(X)$, since otherwise $X$ would not be a minimal endomorphic image of $Z$, a contradiction since $X$ is a core of $Z$. Since $\varphi$ is a homomorphism it preserves adjacency, so we must only check that it preserves non-adjacency on $V(X)$. Suppose not and let $Y$ be the subgraph induced by $\varphi(V(X))$. Any endomorphism of $Z$ must act injectively on $V(Y)$, otherwise composing this endomorphism with $\varphi$ would yield an endomorphism that does not act on $V(X)$ injectively, a contradiction to the above. But then no endomorphism of $Z$ can map $Y$ to $X$, since $Y$ is isomorphic to a graph obtained from $X$ by adding at least one edge. Therefore $\varphi$ must preserve non-adjacency on $V(X)$ and therefore acts as an isomorphism from $X$ to $Y$.\qedn

The following result, proved in Hahn \& Tardif \cite{tardif}, shows that the distance between two vertices in the core of a graph is equal to their distance in the graph. We use $\dist_X(u,v)$ to denoted the \emph{distance in $X$} (length of a shortest path in $X$) between vertices $u$ and $v$. 

\begin{lemma}\label{lem:dist1}
Suppose that $Z$ is a graph with an induced subgraph $X$ isomorphic to the core of $Z$. If $u,v \in V(X)$, then $\dist_X(u,v) = \dist_Z(u,v)$.
\end{lemma}
%\proof
%Let $\rho$ be a retraction of $Z$ onto $X$. Obviously, $\dist_X(u,v) \ge \dist_Z(u,v)$ since $X$ is a subgraph of $Z$. However, the adjacency preserving property of homomorphisms means that they map walks to walks (necessarily of the same length). Thus a shortest path from $u$ to $v$ in $Z$ is mapped by $\rho$ to a walk of equal or lesser length in $X$. Thus $\dist_X(u,v) \le \dist_Z(u,v)$.\qedn

Combining the two above lemmas, we arrive at the following:

\begin{lemma}\label{lem:dist2}
Suppose that $Z$ is a graph with an induced subgraph $X$ isomorphic to the core of $Z$. If $u,v \in V(X)$ and $\varphi$ is an endomorphism of $Z$, then $\dist_Z(u,v) = \dist_Z(\varphi(u), \varphi(v))$.
\end{lemma}
\proof
By Lemma~\ref{lem:iso}, the endomorphism $\varphi$ maps $X$ isomorphically to some other copy, $Y$, of the core of $Z$. Therefore we have $\dist_X(u,v) = \dist_Y(\varphi(u), \varphi(v))$. Together with Lemma~\ref{lem:dist1} we have
\[\dist_Z(u,v) = \dist_X(u,v) = \dist_Y(\varphi(u), \varphi(v)) = \dist_Z(\varphi(u), \varphi(v)).\]  \qedn

As a corollary we have the following:

\begin{corollary}\label{cor:dist1}
Suppose that $\rho$ is a retraction of $Z$ onto a copy $X$ of its core. If $u,v \in V(X)$ and there exists $\sigma \in \aut(Z)$ such that $\sigma(u) = u'$ and $\sigma(v) = v'$, then $\dist_X(\rho(u'),\rho(v')) = \dist_X(u,v)$.
\end{corollary}
\proof
\[\dist_X(\rho(u'),\rho(v')) = \dist_X(\rho\circ \sigma(u), \rho\circ \sigma(v)) = \dist_X(u,v)\]
by Lemmas~\ref{lem:dist1} and~\ref{lem:dist2}, and the fact that $\rho \circ \sigma$ is an endomorphism of $Z$.\qedn

Finally, we can apply the above to the special case of cubelike graphs:

\begin{corollary}\label{cor:dist2}
Let $Z$ be a Cayley graph for $\Ztn$. Let $\rho$ be a retraction of $Z$ onto a copy $X$ of its core. If $a,b \in V(X)$ and $c+d = a+b$, then
\[\dist_X(\rho(c), \rho(d)) = \dist_X(a,b) = \dist_Z(a,b).\]
\end{corollary}
\proof
Adding $a+c$ to all of the elements of $\Ztn$ is an automorphism of $Z$ mapping $a$ to $c$ and $b$ to $d$. Apply above lemmas/corollaries.\qedn

\section{Cubelike Hulls}\label{sec:cubelikehull}

Given a graph $X$, its \emph{cubelike hull}, which is denoted $\ZZ{X}$, is a graph that has the even weight vectors of $\mathbb{Z}_2^{V(X)}$ as its vertex set, and where two vertices are adjacent if their sum is $e_u + e_v$ for some edge $uv$ in $X$.  These even weight vectors form a (sub)group $Z'$, and the cubelike hull is the Cayley graph ${\rm Cay}(Z', \{e_u+e_v: uv \in E(X)\})$. As $Z'$ is isomorphic to $\mathbb{Z}_2^{|V(X)|-1}$, it follows that $\ZZ{X}$ is itself a cubelike graph. 
Letting $C_n$ denote the cycle of length $n$, and $K_n$ the complete graph of order $n$, it is straightforward to confirm that $\ZZ{C_n}$ is the folded cube of order $n$ and $\ZZ{K_n}$ is the halved cube of order $n$.
Note that $X \to \ZZ{X}$ always holds since for any $u \in V(X)$ the vertices $\{e_u+e_v: v \in V(X)\}$ induce a copy of $X$ in $\ZZ{X}$. These graphs were considered by Beaudou, Naserasr, and Tardif in~\cite{BNT15} (where they called them \emph{power graphs}), and they proved the following result, which says that $\ZZ{X}$ is the minimal (in the homomorphism order) cubelike graph admitting a homomorphism from $X$:

\begin{lemma}[Beaudou, Naserasr, and Tardif \cite{BNT15}]\label{lem:Z2}
Let $X$ be a graph and $Z$ a cubelike graph. Then $X \to Z$ if and only if $\ZZ{X} \to Z$. Moreover, the map from $\ZZ{X}$ to $Z$ can be taken to be a group homomorphism.
\end{lemma}
\proof
First, $X \to \ZZ{X}$ always holds and so $\ZZ{X} \to Z$ implies $X \to Z$. Conversely, suppose that $\varphi$ is a homomorphism from $X$ to $Z$. Since $Z$ is cubelike, $V(Z) = \mathbb{Z}_2^d$ for some $d$. Define $\hat{\varphi}: \mathbb{Z}_2^{V(X)} \to V(Z)$ as
\[\hat{\varphi}\left(\sum_{u \in S} e_u\right) = \sum_{u \in S} \varphi(u),\]
where $S$ is an arbitrary subset of $V(X)$. In other words, $\hat{\varphi}$ is the linear extension of the map taking $e_u$ to $\varphi(u)$, and is therefore a group homomorphism. We claim that the restriction of this map to $V(\ZZ{X})$ is a homomorphism from $\ZZ{X}$ to $Z$. Indeed, if $xy \in E(\ZZ{X})$, then $x+y = e_u+e_v$ for some $uv \in E(X)$, and since $\hat{\varphi}$ is linear we have that
\[\hat{\varphi}(x) + \hat{\varphi}(y) = \hat{\varphi}(x+y) = \hat{\varphi}(e_u+e_v) = \hat{\varphi}(e_u)+\hat{\varphi}(e_v) = \varphi(u) + \varphi(v).\]
Since $\varphi$ was a homomorphism and $uv \in E(X)$, we have that $\varphi(u)\varphi(v) \in E(Z)$ which means that $\varphi(u) + \varphi(v)$ is in the connection set of $Z$. Thus $\hat{\varphi}(x)$ and $\hat{\varphi}(y)$ are adjacent as desired.\qedn

The primary focus of this work is the question of whether the core of a cubelike graph is necessarily cubelike. Given a cubelike graph $Z$, it is straightforward how to determine whether its core is cubelike: first find its core $Z^\bullet$, then decide whether it is cubelike by testing whether its automorphism group contains a regular elementary abelian 2-group. However, it is not so clear how to go in the other direction: given a core $X$, how does one determine whether it is the core of \emph{some} cubelike graph? Indeed, it is not even clear if this problem should be decidable, since in principle one may have to consider every cubelike graph $Z$ and determine whether $Z^\bullet \cong X$. However, the following two lemmas, based on Lemma~\ref{lem:Z2}, provide us with an algorithm for this problem\footnote{After arriving at these results independently, we learned that they appear in~\cite{homotens}.}.

\begin{lemma}
A graph $X$ is homomorphically equivalent to a cubelike graph if and only if it is homomorphically equivalent to $\mathbb{Z}_2[X]$.
\end{lemma}
\proof
If $X$ is homomorphically equivalent to $\mathbb{Z}_2[X]$, then it is homomorphically equivalent to a cubelike graph since the latter is cubelike.

Conversely, suppose that $X$ is homomorphically equivalent to a cubelike graph $Z$. Since $X \to Z$ and $Z$ is cubelike, we have that $\mathbb{Z}_2[X] \to Z \to X$ by Lemma~\ref{lem:Z2}. As noted above, $X \to \mathbb{Z}_2[X]$ is true for any $X$ and so $\mathbb{Z}_2[X]$ and $X$ are homomorphically equivalent.\qedn

\begin{lemma}\label{lem:Z2core}
A graph $X$ is the core of a cubelike graph if and only if it is the core of $\mathbb{Z}_2[X]$.
\end{lemma}
\proof
If $X$ is the core of $\mathbb{Z}_2[X]$, then it is the core of a cubelike graph since the latter is cubelike.

Conversely, if $X$ is the core of a cubelike graph, then by the above lemma $X$ must be homomorphically equivalent to $\mathbb{Z}_2[X]$. But since $X$ is a core, it must be the core of $\mathbb{Z}_2[X]$.\qedn

By the above, in order to determine whether a particular graph $X$ is the core of \emph{some} cubelike graph, we can check whether it is a core and whether $\ZZ{X} \to X$. Unfortunately, this is not feasible in practice because if $X$ is any graph with more than $8$ vertices that might be the core of a cubelike graph, then $X$ has at least $16$ vertices and  $\ZZ{X}$ has at least $2^{15} = 32768$ vertices. However, our results will allow us to prove by hand that if $X$ is the core of a cubelike graph and $|V(X)| \le 16$, then $X$ is itself cubelike. Additionally, with the aid of a computer we are able to prove this for $|V(X)| \le 32$.

On the other hand, Lemma~\ref{lem:Z2} does provide us with useful necessary conditions for the core of a cubelike graph, which we will use to prove these results. Indeed we have the following:

\begin{lemma}\label{lem:Z2ext}
If $X$ is the core of a cubelike graph, then $Y \to X$ if and only if $\ZZ{Y} \to X$
\end{lemma}
\proof
Suppose that $X$ is the core of a cubelike graph $Z$. If $\ZZ{Y} \to X$, then $Y \to X$ since $Y \to \ZZ{Y}$. Conversely, suppose that $Y \to X$. Since $X$ and $Z$ are homomorphically equivalent, we have that $Y \to Z$ and thus by Lemma~\ref{lem:Z2} we have that $\ZZ{Y} \to Z$ and thus $\ZZ{Y} \to X$, again using the homomorphic equivalence of $X$ and $Z$.\qedn

This lemma will be of significant use to us in Section~\ref{sec:fewverts}. In particular, we will apply Lemma~\ref{lem:Z2ext} when $Y$ is a complete graph (so $\ZZ{Y}$ is a halved cube) or when $Y$ is an odd cycle (so $\ZZ{Y}$ is a folded cube). The next lemma extends results of Beaudou, Naserasr and Tardif~\cite{BNT15}.

%that $\ZZ{C_n}$ is isomorphic to the folded cube of order $n$ (there they refer to it as the \emph{projective cube of dimension $n$}). There they also prove that if a cubelike graph has odd girth $n$, then it must contain the folded cube of order $n$ as a subgraph. It is simple to additionally prove that it must contain this folded cube as an \emph{induced} subgraph, and that this holds even for graphs that are homomorphically equivalent to a cubelike graph. We will use this extended version of their result, and provide a proof as well.

%\begin{lemma}
%For $g \ge 2$, the graph $\mathbb{Z}_2[C_{g}]$ is isomorphic to the folded cube of order $g$.
%\end{lemma}

\begin{lemma}\label{lem:FQsub}
Let $X$ be a graph that is homomorphically equivalent to a cubelike graph $Z$. If $X$ has odd girth $g < \infty$, then it contains $\mathbb{Z}_2[C_{g}]$ as an induced subgraph.
\end{lemma}
\proof
Since $X$ has odd girth $g$, we have that $C_g \to X$. Since $X$ is homomorphically equivalent to a cubelike graph, we have that $\mathbb{Z}_2[C_{g}] \to X$. Let $\varphi$ be a homomorphism from $\mathbb{Z}_2[C_{g}]$ to $X$. Since any two vertices of $\ZZ{C_g}$ are contained in a cycle of length $g$, and identifying two vertices of an odd cycle results in a graph containing an odd cycle of strictly shorter length, the map $\varphi$ must be injective, i.e., $\ZZ{C_g}$ is a subgraph of $X$. Furthermore, adding an edge between any pair of non-adjacent vertices in an odd cycle also results in a graph containing a shorter odd cycle. Therefore, the map $\varphi$ must be faithful, i.e., $\ZZ{C_g}$ must be an induced subgraph of $X$.\qedn

%We will argue that $\varphi$ is injective and faithful.
%
%Suppose $u,v$ are distinct vertices of $\mathbb{Z}_2[C_{g}]$ such that $\varphi(u) = \varphi(v)$. Since any two vertices of a folded cube are contained in a shortest odd cycle, we have that $u$ and $v$ are contained in some cycle $C$ of length $g$ in $\mathbb{Z}_2[C_{g}]$. The restriction of $\varphi$ to $C$ is still a homomorphism to $X$, and since $\varphi(u) = \varphi(v)$ this homomorphism is not injective. This implies that the image of $C$ under $\varphi$ must contain a shorter odd cycle which contradicts the fact that $X$ has odd girth $g$. Thus $\varphi$ must be injective.
%
%Suppose that $\varphi$ is not faithful, i.e., there exist adjacent vertices $x$ and $y$ of $X$ that are in the image of $\varphi$ but there are no edges between $\varphi^{-1}(x)$ and $\varphi^{-1}(y)$. Let $u \in \varphi^{-1}(x)$ and $v \in \varphi^{-1}(y)$ and let $C$ be a shortest odd cycle of $\mathbb{Z}_2[C_{g}]$ containing $u$ and $v$. The image of $C$ under $\varphi$ is therefore an odd cycle of length $g$ containing $x$ and $y$, but $x$ and $y$ are not adjacent in this cycle. Thus, the edge $xy$ is a chord of $\varphi(C)$, and it follows that $X$ contains a shorter odd cycle, a contradiction.\qedn

The above lemma is useful since looking for induced subgraphs isomorphic to a graph $X$ is often easier than looking for homomorphisms from $X$. This will also be of theoretical use, for instance in the proof of Theorem~\ref{thm:fewverts}.

\section{Shift Graphs and Hom-Idempotence}

The \emph{Cartesian product} of graphs $X$ and $Y$, denoted $X \cart Y$, has vertex set $V(X) \times V(Y)$, and there is an edge between $(x,y)$ and $(x',y')$ if $x = x'$ and $y \sim y'$, or $x \sim x'$ and $y = y'$. Note that $X \to X \cart Y$ since for any fixed $y \in V(Y)$ the subset $\{(x,y) : x \in V(X)\}$ induces a copy of $X$ in $X \cart Y$, and similarly $Y \to X \cart Y$. Therefore, for any graph $X$, we have that $X \to X \cart X$ and a natural question to ask is whether/when $X \cart X \to X$ also holds. This question was considered by Larose, Laviolette, and Tardif in~\cite{LLT}, where they called graphs $X$ with the property that $X \cart X \to X$, \emph{hom-idempotent}. They noted that if $X = \text{Cay}(\Gamma,C)$ is a normal Cayley graph, then the map $(g,h) \mapsto gh$ is a homomorphism from $X \cart X$ to $X$. Thus normal Cayley graphs are always hom-idempotent. Moreover, since a homomorphism $\varphi$ from $X$ to $Y$ can be used to construct a homomorphism $(x,x') \mapsto (\varphi(x),\varphi(x'))$ from $X \cart X$ to $Y \cart Y$, as long as $X$ is homomorphically equivalent to a normal Cayley graph it will be hom-idempotent. Remarkably, Larose et.~al.~prove the converse: that a graph $X$ is hom-idempotent if and only if it is homomorphically equivalent to a normal Cayley graph. Moreover, they showed that $X$ is hom-idempotent if and only if it is homomorphically equivalent to a \emph{particular} normal Cayley graph.

\begin{definition}
Given a graph $X$, an automorphism $\sigma \in \aut(X)$ is a \emph{shift} of $X$ if $\sigma(x) \sim x$ for all $x \in V(X)$, i.e., $\sigma$ maps each vertex to one of its own neighbors. The \emph{shift graph} of $X$, denoted $\text{Sh}(X)$, is the Cayley graph $\text{Cay}(\aut(X),S)$ where $S$ is the set of shifts of $X$.
\end{definition}

Note that the inverse of a shift is a shift ($x \sim \sigma(x) \Rightarrow \sigma^{-1}(x) \sim \sigma^{-1}(\sigma(x)) = x$), and therefore the shift graph of $X$ is indeed a Cayley \emph{graph}. Moreover, if $\sigma$ is a shift and $\pi$ is an automorphism, then $\pi \circ \sigma \circ \pi^{-1} (\pi(x)) = \pi(\sigma (x)) \sim \pi(x)$ since $\sigma(x) \sim x$ and $\pi$ is an automorphism. By varying $x$ over all of $V(X)$, we also vary $\pi(x)$ over all of $V(X)$, and thus $\pi \circ \sigma \circ \pi^{-1}$ is a shift. Thus $\text{Sh}(X)$ is a normal Cayley graph. Larose et.~al.~proved the following~\cite{LLT}:

\begin{lemma}\label{lem:homidem}
Let $X$ be a graph. Then the following are equivalent:
\begin{enumerate}
\itemsep0em 
\item $X$ is hom-idempotent;
\item $X$ is homomorphically equivalent to a normal Cayley graph;
\item $X$ is homomorphically equivalent to $\mathrm{Sh}(X^\bullet)$.
\end{enumerate}
\end{lemma}

For any graph $X$ and vertex $x \in V(X)$, the map $\sigma \mapsto \sigma^{-1}(x)$ is a homomorphism from $\text{Sh}(X)$ to $X$. Therefore, if $X$ is a core then it is homomorphically equivalent to its shift graph if and only if the latter contains the former as an induced subgraph.

The above lemma provides two algorithms for checking whether a given graph $X$ is homomorphically equivalent to some normal Cayley graph: either check if $X$ is hom-idempotent, or check whether $X^\bullet$ is a subgraph of $\text{Sh}(X^\bullet)$. This is noteworthy since it is not a priori obvious that such an algorithm should exist, as with the question of whether a graph is homomorphically equivalent to some cubelike graph, which was discussed in Section~\ref{sec:cubelikehull}.

It is not at all obvious from the definition that hom-idempotence should have anything to do with normal Cayley graphs. Thus part of the significance of the above lemma is that it hints at a deep connection between homomorphisms and normal Cayley graphs, thus providing some justification for Question~\ref{qn:normcore}.

\section{Orbital Graphs}

Given a permutation group $\Gamma$ acting on a set $V$, the \emph{orbit} of an element $v \in V$ is $\Gamma(v) = \{\gamma(v) : \gamma \in \Gamma\}$. The orbits of $\Gamma$ partition $V$ and, almost by definition, $\Gamma$ is transitive if and only if it has a single orbit. The action of $\Gamma$ on $V$ also induces a natural action on $V \times V$, and the orbits of this action are referred to as the \emph{orbitals} of $\Gamma$, and these partition the set $V \times V$. 

If $O$ is an orbital, then the set  
$
O^* = \{(y,x) : (x,y) \in O\},
$
is also an orbital, called the \emph{paired} orbital of $O$. If $O^*= O$ then $O$ is a \emph{self-paired} orbital. Any orbital that contains only pairs of the form $(x,x)$ is called a \emph{diagonal} orbital; there is a unique diagonal orbital if and only if  $\Gamma$ is a transitive group.

As with orbits, certain properties of the group $\Gamma$ can be characterized in terms of its orbitals. For instance, a transitive group is generously transitive if and only if all of its orbitals are self-paired, and it is regular if and only if it has exactly $|X|$ orbitals.  Each orbital can be viewed as the set of directed edges of a directed graph, and even more properties of a transitive group can be determined by considering the digraphs whose edges are the unions of orbitals.

\begin{definition}
Let $\Gamma$ be a permutation group acting on a set $V$.  An \emph{orbital digraph} of $\Gamma$ is a digraph with vertex set $V$ and whose arc set $A$ is the union of non-diagonal orbitals of $\Gamma$. If $X$ is a (di)graph, then we will additionally refer to the orbital digraphs of $\aut(X)$ as the orbital digraphs of $X$.
\end{definition}

%Note that the above definition allows for orbital (di)graphs to have loops. Since any graph admits a homomorphism to any graph with a loop, this case will not be particularly interesting for us, though our results will still hold. 

If the set of orbitals whose union is $A$ is closed under taking paired orbitals, the resulting orbital digraph will contain a given directed edge if and only if it contains the reverse edge. In this case, we can consider them as (undirected) graphs and will refer to them simply as \emph{orbital graphs}. Since we are mostly concerned with cubelike graphs, which are generously transitive and have only self-paired orbitals, this will be the usual case for us.   

%We remark that a transitive permutation group is primitive if and only if all of its non-diagonal orbital digraphs are (strongly) connected.

The following theorem shows that there is a connection between the orbital digraphs of a graph $X$ and those of its core $X^\bullet$. Below, we use $\Gamma(E)$ to denote the set $\{(\gamma(x),\gamma(y)) : (x,y) \in E, \ \gamma \in \Gamma\}$, where $E$ is a set of ordered pairs of elements of the set $\Gamma$ acts on. Recall that for $S \subseteq V(X) \times V(X)$ the notation $X(S)$ refers to the (di)graph with vertex set $V(X)$ and arc/edge set $S$.

\begin{theorem}\label{thm:orbits}
Let $X$ be a graph and let $\rho$ be a retraction of $X$ onto a copy $X^\bullet$ of its core. Suppose that $E$ is a union of orbitals of $\aut(X^\bullet)$ and that $\Gamma \le \aut(X)$. Then $\rho$ is a homomorphism from $X(\Gamma(E))$ to the digraph $X^\bullet(E)$. In particular, $X(\Gamma(E))$ and $X^\bullet(E)$ are homomorphically equivalent.
\end{theorem}

%Let $Z$ be a Cayley graph for $\Ztn$ and let $X$ be a copy of its core. Let $O_1, \ldots, O_m$ be the orbits of the automorphism group of $X$ on the pairs of vertices of the core (the automorphism group is generously transitive and so we do not need to consider ordered pairs of vertices). For each $i \in \{1, \ldots, m\}$ let $D_i = \{a+b : \{a,b\} \in O_i\}$, i.e.~$D_i$ is the set of differences (sums) of the pairs of vertices appearing in orbit $O_i$. For any $S \subseteq [m]$ let $Z_S$ be the Cayley graph on $\Ztn$ with connection set $\cup_{i \in S} D_i$. Then $Z_S$ is homomorphically equivalent to the graph with vertex set $V(X)$ and edge set $\cup_{i \in S} O_i$.
\proof
Suppose that $xy$ is an arc of $X(\Gamma(E))$. We must show that $\rho(x)\rho(y)$ is an arc of $X^\bullet(E)$. By the definition of $X(\Gamma(E))$, there exists $ab \in E$ and $\gamma \in \Gamma$ such that $\gamma(a) = x$ and $\gamma(b) = y$. Since $\gamma$ is an automorphism of $X$, the map $\rho \circ \gamma$ is an endomorphism of $X$ whose image is $X^\bullet$. Thus, the restriction $\rho \circ \gamma \restriction_{V(X^\bullet)}$ is an automorphism of $X^\bullet$, and therefore $\rho\circ \gamma (ab) \in E$ since $ab \in E$. This implies that $\rho(x)\rho(y)$ is an arc of $X^\bullet(E)$, and so we have shown that $X(\Gamma(E)) \to X^\bullet(E)$. We trivially have $X^\bullet(E) \to X(\Gamma(E))$ since the former is a subgraph of the latter.\qedn

If we apply the above theorem to Cayley graphs we obtain the following:

\begin{corollary}\label{cor:cayorbits}
Let $X$ be a Cayley (di)graph for a group $\Gamma$. Then every orbital (di)graph of $X^\bullet$ is homomorphically equivalent to a Cayley (di)graph for $\Gamma$.
\end{corollary}
\proof
Since $X$ is a Cayley graph for $\Gamma$, we have that $\Gamma \le \aut(X)$ and $\Gamma$ acts regularly on $X$. If $E$ is any union of orbitals of $\aut(X^\bullet)$, then taking $\Gamma$ in the above theorem to be this regular subgraph results in a graph $X(\Gamma(E))$ with $\Gamma$ acting regularly on its vertices, i.e., it is a Cayley graph for $\Gamma$.\qedn

This is rather interesting, as well as useful for our work. For example, it implies that if $X$ is the core of a cubelike graph, then no orbital graph of $X$ can have clique number equal to three. Indeed, any orbital graph of $X$ is homomorphically equivalent to some cubelike graph $Z$ by the above corollary, and thus $\omega(X) = \omega(Z) \ne 3$. In practice, this a very useful tool for ruling out potential counterexamples to Conjecture~\ref{conj:cubecore}, as we will see in Section~\ref{subsec:32vertices}.

\section{The Covering Cube Theorem and the Degree Bound}

In this section we will show that if $Z$ is a cubelike graph with core $X$, then there is a cubelike subgraph of $Z$ that covers $X$. Moreover, this implies that $X$ is covered by a $d$-cube where $d$ is the degree of $X$. This latter statement follows from the former due to the fact that any connected cubelike graph of degree $d$ is covered by the $d$-cube, which we prove in Lemma~\ref{lem:cubecover} below. Being covered by the $d$-cube further implies a useful necessary condition for a graph to be the core of a cubelike graph: its eigenvalues must be a sub-multiset of the eigenvalues of the $d$-cube (Corollary~\ref{cor:spectra}). In fact, this applies to all of the orbital graphs of the core of a cubelike graph, making this requirement even more stringent.

\begin{lemma}\label{lem:cubecover}	
Let $Z$ be a connected cubelike graph with degree $d$. Then $Z$ is covered by the $d$-cube.
\end{lemma}
\proof
By assumption $Z = \text{Cay}(\Ztn,C)$ for some $n \in \mathbb{N}$ and $C \subseteq \Ztn \setminus \{0\}$ such that $|C| = d$. Let $C = \{c_1, \ldots, c_d\}$. Recall that the $d$-cube is the graph $Q_d = \text{Cay}(\mathbb{Z}_2^d,\{e_1, \ldots, e_d\})$. Considering $\Ztn$ and $\mathbb{Z}_2^d$ as vector spaces, let $f : \mathbb{Z}_2^d \to \Ztn$ be the linear extension of the map defined as $f(e_i) = c_i$. We will show that $f$ is a covering map from $Q_d$ to $Z$.

First, since $Z$ is connected, we have that the set $\{c_1, \ldots, c_d\}$ spans all of $\Ztn$. This implies that $f$ is surjective. Since $Z$ has the same degree as $Q_d$, it remains only to show that $f$ is locally injective. Let $x$ be a vertex of $Q_d$, and consider two of its neighbors $x+e_i$ and $x+e_j$ for $i \ne j$. We have that $f(x+e_i) = f(x) + f(e_i) = f(x) + c_i$, and similarly $f(x+e_j) = f(x) + c_j$. Note that this shows that $f$ maps the neighbors of $x$ to the neighbors of $f(x)$, i.e., $f$ is a homomorphism. Furthermore, $f(x+e_i) + f(x + e_j) = c_i + c_j \ne 0$ since $i \ne j$. This implies that $f(x + e_i) \ne f(x + e_j)$ and so $f$ is locally injective as desired.\qedn

In Theorem~\ref{thm:cover} below, we will further show that the core of a cubelike graph is covered by the cube (of the appropriate degree). The motivation behind this result came from considering how one might identify all the graphs of a particular fixed degree that might be the core of a cubelike graph.  As an illustrative example, suppose that we would like to determine all of the graphs  of degree $3$ that are cores of some cubelike graph. How might we go about this?\footnote{For $d=3$, we can alternatively use Brooks' theorem which says that such a graph $X$ either has chromatic number at most 3 or is $K_4$. If the chromatic number is 2, then $X$ is bipartite and thus not a core unless it is $K_2$, which does not have degree 3. Thus the chromatic number of $X$ must be 3, but no cubelike graph can have chromatic number 3 by the result of Payan~\cite{payan}, and thus neither can any core of a cubelike graph. Therefore the only cubic core of a cubelike graph is $K_4$. This is quicker than the method described, but does not generalize to higher valencies.}

Suppose $Z$ is a cubelike graph on $\mathbb{Z}_2^n$ and $\varphi$ is a retraction onto a copy $X$ of its core which contains the zero element. Further suppose that $X$ has degree 3 and that the neighbors of $0$ in $X$ are $a_1,a_2,a_3 \in \mathbb{Z}_2^n$. Note that this implies that $a_1,a_2,a_3$ are in the connection set of $Z$. Since $\varphi$ is a retraction, we have that $\varphi(0) = 0$ and $\varphi(a_i) = a_i$ for $i = 1,2,3$. Consider $\varphi(a_1+a_2)$. Since $a_1 + (a_1+a_2) = a_2$, we have that $a_1+a_2$ is adjacent to $a_1$, and similarly $a_2$. Thus $\varphi(a_1+a_2)$ must be adjacent to $\varphi(a_1) = a_1$ and $\varphi(a_2) = a_2$. Also, note that $a_1$ and $a_2$ are at distance 1 or 2 in the core, depending on whether they are adjacent. Since $0 + (a_1+a_2) = a_1+ a_2$, Corollary~\ref{cor:dist2} implies that $\varphi(a_1+a_2)$ is at distance 1 or 2 from $\varphi(0) = 0$. In particular, $\varphi(a_1+a_2) \ne 0$. There are thus two possibilities for $\varphi(a_1+a_2)$, either it is not among $0,a_1,a_2,a_3$, or it is equal to $a_3$. 

Suppose that $\varphi(a_1+a_2) = a_3$. We will show that $a_i+a_j$ is in the connection set for all $i \ne j \in \{1,2,3\}$. Since $\varphi(a_1+a_2) = a_3$, we have that $a_1$ and $a_2$ must be adjacent to $a_3$ since they are adjacent to $\varphi(a_1+a_2)$. This implies that $a_1+a_3$ and $a_2+a_3$ are in the connection set of $Z$, and thus $0$ is adjacent to $a_2+a_3$ and $\varphi(a_2+a_3)$.  Therefore $\varphi(a_2+a_3) \in \{a_1,a_2,a_3\}$, and it follows that $\varphi(a_2+a_3) = a_1$ since the other two vertices must be adjacent to $a_2+a_3$. But this now implies that $a_1 = \varphi(a_2+a_3)$ and $a_2$ are adjacent, and so $a_1+a_2$ is in the connection set of $Z$. Thus we have shown that $a_i+a_j$ is in the connection set for all $i \ne j$. This implies that the vertices $0,a_1,a_2,a_3$ induce a $K_4$ in the core $X$. However, $X$ has degree 3 by assumption and must be connected since it is a vertex-transitive core, and so this $K_4$ must be all of $X$. Note that $K_4$ is indeed the core of a cubelike graph, for instance $K_4$ itself is cubelike.

Now suppose that $\varphi(a_1+a_2) \ne a_3$. By symmetry we may assume that none of $\varphi(a_i+a_j)$ for $i \ne j$ are equal to any of $0,a_1,a_2,a_3$, since otherwise we would be in the previous case. Note that this implies that none of $a_1,a_2,a_3$ are adjacent, since if $a_1$ and $a_2$ were adjacent, for instance, then $a_1+a_2$ would be in the connection set, and so $0$ would be adjacent to $a_1+a_2$. This would imply that 0 is adjacent to $\varphi(a_1+a_2)$ and thus the latter vertex must be among the neighbors of $0$ in $X$, i.e., must be one of $a_1,a_2,a_3$. Therefore, the distance between $a_i$ and $a_j$ in $X$ is exactly 2 for $i \ne j$, and so by Corollary~\ref{cor:dist2} we have that $\varphi(a_i+a_j)$ is at distance 2 from $\varphi(a_\ell+a_k)$ whenever $\{i,j\} \ne \{\ell,k\}$. In particular this implies that the three vertices of the form $\varphi(a_i+a_j)$ for $i \ne j$ are distinct. Also note that since $(a_i+a_j) + a_i = a_j$ which is in the connection set, we have that $\varphi(a_i+a_j)$ is adjacent to both $a_i$ and $a_j$. Finally, the vertex $\varphi(a_1+a_2+a_3)$ must be adjacent to the three vertices of the form $\varphi(a_i+a_j)$. Considering the adjacencies we have shown among the eight vertices we have considered so far, it is not difficult to see that these form a 3-cube. Since the 3-cube is regular of degree 3, and the core $X$ is regular of degree 3 and connected, this must be all of $X$. However, this is a contradiction since the 3-cube is bipartite and thus not a core.

So we have shown that the only degree 3 core of a cubelike graph is $K_4$. The above approach can also be carried through for degree 4 or even 5, though it becomes quite tedious at this point. The truly determined may even be able to do the degree 6 or 7 case, though we do not recommend attempting this. The idea of the above is that we consider the subgroup $\Gamma$ of $\mathbb{Z}_2^n$ generated by the neighbors of $0$ in the core, and then use Corollary~\ref{cor:dist2} and case analysis to reason about the possible images of the elements of $\Gamma$. One might wonder whether all of the core is necessarily revealed once one determines the images of $\Gamma$ and the adjacencies among them. After some thought one comes to the following theorem:

\begin{theorem}\label{thm:cover}
Let $Z$ be a Cayley graph for $\Ztn$ and let $\rho$ be a retraction onto a core, $X^\bullet$, which contains $0$. Let $X$ be an orbital graph of $X^\bullet$ and let $a_1, \ldots, a_d$ be the neighbors of 0 in $X$. Let $Y = \mathrm{Cay}(\langle a_1, \ldots, a_d \rangle, \{a_1, \ldots, a_d\})$. Then the restriction, $\rho\restriction_Y$, is a covering map from $Y$ to the connected component of $X$ containing 0. Moreover, this implies that this connected component is covered by the $d$-cube.
\end{theorem}
\proof
Let $O$ be the union of orbitals of $X^\bullet$ such that $X = X^\bullet(O)$, and let $X'$ be the connected component of $X$ containing the vertex 0. Since $\Ztn \le \aut(Z)$ and $(0,a_i) \in O$ for $i \in [d]$, we have that $Y$ is a subgraph of the orbital graph of $Z$ with edge set $\Ztn(O)$. Therefore, by Theorem~\ref{thm:orbits} the map $\rho\restriction_Y$ is a homomorphism to $X$. However, $Y$ is connected and thus its image under $\rho$ must be connected. Since $Y$ contains the vertex 0 and $\rho$ is a retraction, we see that $\rho\restriction_Y$ is a homomorphism to $X'$. It remains to show that $\rho\restriction_{Y}$ is surjective and locally bijective. We first prove the latter.

Since $Y$ has the same degree as $X'$ by construction, it suffices to show that it is locally injective. Consider $y \in V(Y)$ and its neighbors $y+a_1, \ldots, y+a_d$ in $Y$. For distinct $i,j \in [d]$, we have that $(y+a_i) + (y+a_j) = a_i + a_j$. Since $a_i$ and $a_j$ appear in $X^\bullet$ at nonzero distance from each other, we have that $\rho(y+a_i) \ne \rho(y+a_j)$ by Corollary~\ref{cor:dist2}. Therefore $\rho\restriction_Y$ is locally injective.

Since $\rho\restriction_Y$ is locally injective, the image of $Y$ under this map is a subgraph of $X'$ of minimum degree at least $d$. However, since $X'$ is connected, this implies that the image must be all of $X'$. Therefore, the map $\rho\restriction_Y$ is surjective and thus we have shown it is a covering map.

Finally, by Lemma~\ref{lem:cubecover}, the graph $Y$ is covered by the $d$-cube. Since the composition of covering maps is a covering map, we have that $X'$ is covered by the $d$-cube.\qedn

The covering cube theorem has two immediate consequences. The first gives a bound on the number of vertices in the core of a cubelike graph in terms of the degree of the core. We refer to this as the \emph{degree bound}:

\begin{corollary}\label{cor:degreebound}
Let $X$ be the core of a cubelike graph. If $X'$ is a connected component of an orbital graph of $X$ with degree $d$, then $X'$ has at most $2^d$ vertices. If $X$ has degree $d$, then it has at most $2^{d-1}$ vertices unless $d = 1$ and $X = K_2$.
\end{corollary}
\proof
Since $X$ is the core of a cubelike graph, it is vertex transitive and thus so are all of its orbital graphs. Thus the connected components of any fixed orbital graph are all isomorphic. If an orbital graph $X'$ of $X$ has degree $d$, then each of its connected components are covered by the $d$-cube by Theorem~\ref{thm:cover} above. The $d$-cube has $2^d$ vertices, and thus each component of $X'$ has at most this many vertices.

In the case of $X$ itself, if it has $2^d$ vertices then the covering map from the $d$-cube must in fact be an isomorphism. However, this implies that $X$ is bipartite and therefore not a core unless $X = K_2$, in which case $d = 1$. Therefore, $X$ has at most $2^{d-1}$ vertices unless $d = 1$ and $X = K_2$.\qedn

Note that if we apply the above to the case of degree 3 cores of cubelike graphs, we immediately obtain that such a core has at most 4 vertices. Since it has degree 3 it must also have at least 4 vertices and it must be $K_4$. So we get that the only degree 3 core of a cubelike graph is $K_4$, and this was much quicker than the argument given before Theorem~\ref{thm:cover}. Moreover, this argument avoids using Payan's nontrivial result about the chromatic number of cubelike graphs (which we used in the alternative argument given in the footnote).

Recall from Section~\ref{sec:preliminaries} that if a graph $X$ covers a graph $Y$, then the eigenvalues of $Y$ are a submultiset of the eigenvalues of $X$. Therefore, by the Covering Cube Theorem we have the following:

%requires us to introduce \emph{equitable partitions}. A partition $\{P_1, \ldots, P_k\}$ of the vertex set of a graph $X$ is \emph{equitable} if there exist numbers $c_{ij}$ such that for all $i,j \in [k]$, every vertex in $P_i$ has $c_{ij}$ neighbors in $P_j$. The matrix $C = (c_{ij})$ is called the \emph{quotient matrix} of the partition, and it is known~\cite[Lemma~2.3.1]{spectraofgraphs} that the eigenvalues of $C$ are a sub-multiset of the eigenvalues of the adjacency matrix of $X$. From this and the Covering Cube Theorem we obtain the following:

\begin{corollary}\label{cor:spectra}
Let $X$ be the core of a cubelike graph. If $X'$ is a connected component of an orbital graph of $X$ with degree $d$, then its eigenvalues are a sub-multiset of the eigenvalues of the $d$-cube.
\end{corollary}
%\proof
%Let $\varphi$ be the covering map from the $d$-cube to $X'$. We will partition the vertices of the $d$-cube according to the preimages of $\varphi$, i.e., we will consider the partition $\{\varphi^{-1}(x) : x \in V(X')\}$. We will show that this is equitable. Consider two vertices $x,y \in V(X')$. If $x \not\sim y$, then no vertex in $\varphi^{-1}(x)$ is adjacent to any vertex of $\varphi^{-1}(y)$, since $\varphi$ is a homomorphism. On the other hand, if $x \sim y$, then each vertex of $\varphi^{-1}(x)$ is adjacent to exactly one vertex of $\varphi^{-1}(y)$ since $\varphi$ is locally bijective. This implies that our partition is equitable, and its quotient matrix is exactly the adjacency matrix of $X'$. Therefore, the eigenvalues of $X'$ are a sub-multiset of the eigenvalues of the $d$-cube.\qedn

The eigenvalues of the $d$-cube are $d-2i$ for $i = 0, 1, \ldots, d$ with respective multiplicities $\binom{d}{i}$. Thus the core of a cubelike graph, and its orbital graphs, must all have integer eigenvalues. 

We remark that the Covering Cube theorem does not characterize cores of cubelike graphs: there are cores that are covered by a cube but are not the core of any cubelike graph. For instance, let $S$ be the \emph{Shrikhande graph} and $R_{4,4}$ be the $4 \times 4$ \emph{rook graph} (the Cartesian product $K_4 \cart K_4$). Each is a $(16,6,2,2)$ strongly regular graph, hence they have the same spectrum despite being non-isomorphic. The Shrikhande graph is a core, but it has clique number three, so it cannot be the core of any cubelike graph. The $4 \times 4$ rook graph is cubelike (though not a core). Since $R_{4,4}$ is cubelike, so is its \emph{bipartite double cover}, and so they are both covered by the 6-cube. However the bipartite double cover of the Shrikhande graph is isomorphic to that of $R_{4,4}$. Therefore, $S$ is also covered by the 6-cube, even though it is not the core of any cubelike graph.

\subsection{The equality case of the degree bound}

As with many bounds, it is of interest to consider what can be said in the case where the degree bound is met with equality. It is easy to see that if an orbital graph of the core of a cubelike graph has degree $d$ and $2^d$ vertices, then the covering map from the $d$-cube to the orbital graph must be an isomorphism. However, in the case of the core itself, where the bound is $2^{d-1}$ vertices, equality does not imply that the covering map is an isomorphism. In fact, in this case the map will have fibres of size two. We will see that this can still be handled, though with significantly more work.

We will first need the following lemma.

\begin{lemma}\label{lem:cosets}
Let $Y$ be a Cayley graph for $\Ztn$ and $\varphi$ a vertex- and edge-surjective homomorphism from $Y$ to a graph $X$ such that the fibres of $\varphi$ are all cosets of a subgroup $\Gamma$ of $\Ztn$. Then $X$ is a cubelike graph.
\end{lemma}
\proof
Let $X$, $Y$, $\Gamma$, and $\varphi$ be as stated, and let $C$ be the connection set of $Y$. Since the fibres of $\varphi$ are the cosets of $\Gamma$, we can label the vertices of $X$ with these cosets, which form the quotient group $\Ztn/\Gamma$. Define $C' = \{y+\Gamma : \exists g \in \Gamma \text{ s.t.~} y+g \in C\} \subseteq \Ztn/\Gamma$. Note that $0+\Gamma \not\in C'$ since $\varphi$ is a homomorphism and therefore its fibres are independent sets.

Since $\varphi$ is edge surjective, we have that $a+\Gamma \sim b+\Gamma$ in $X$ if and only if there exists $g,g' \in \Gamma$ such that $a+g \sim b+g'$ in $Y$ if and only if $(a+b) + (g + g') \in C$ if and only if $C' \ni (a+b)+\Gamma = (a+\Gamma) + (b+\Gamma)$. Therefore, $X$ is the Cayley graph of $\Ztn/\Gamma$ with connection set $C'$. Since $\Ztn/\Gamma$ must be isomorphic to some power of $\mathbb{Z}_2$, we have proven the lemma.\qedn

To characterize the graphs which meet the degree bound, we first need to show that the only cubelike cores which meet it are the folded cubes of odd order.

\begin{lemma}\label{lem:cubecores}
Let $X$ be a cubelike graph that is a core, has degree $d$, and has $2^{d-1}$ vertices. Then $d$ is odd and $X$ is the folded cube of order $d$.
\end{lemma}
\proof
Since $X$ is a core it is connected and therefore its connection set $C$ contains a basis. Without loss of generality we may assume that $C$ contains the $d-1$ standard basis vectors and one other element $c$. Suppose that $c$ is not the all ones vector. We may assume that $c = \sum_{i=1}^m e_i$ for some $m < d-1$. We will show that this is not a core.

We define a proper endomorphism $\varphi$ of $X$ as follows. Let $\varphi$ be the linear extension of the map fixing $e_i$ for $i \le m$ and sending $e_i$ to $c$ for $i > m$. Note that this maps every element of $C$ to an element of $C$. Suppose that $x \sim y$, i.e., $x+y \in C$. Then $\varphi(x) +\varphi(y) = \varphi(x+y) \in C$, and thus $\varphi(x) \sim \varphi(y)$. Therefore, $\varphi$ is an endomorphism and it clearly is not surjective (nothing is mapped to $e_i$ for $i > m$), a contradiction to the assumption that $X$ is a core.

%The element $\varphi(x)$ is obtained from the element $x$ by changing the last $d-1-m>0$ coordinates to 0 and then adding $c$ if these coordinates contained an odd number of 1's. It is easy to see that if $x,y \in \Ztn$ are such that $x+y = a_i$ for $i \le m$, then $\varphi(x)+\varphi(y) = a_i$, and if $x+y = c$ then $\varphi(x) +\varphi(y) = c$. Moreover, if $x+y = a_i$ for $i > m$, then $\varphi(x)+\varphi(y) = c$. This implies that $\varphi$ is an endomorphism, and it is clearly not surjective. Therefore, $X$ is not a core, a contradiction.

%Then it is easy to see that $X$ is isomorphic to the Cartesian product of the folded cube of order $m+1$ and the $(k-m-1)$-cube. However, this implies that $X$ has a homomorphism to the folded cube of order $m+1$ which it also contains as a subgraph, a contradiction to the fact that $X$ is a core.

The above implies that the connection set of $X$ must consist of the $d-1$ standard basis vectors and the all ones vector, which is one of the standard constructions of the folded cube of order $d$. If $d$ is odd then this graph is indeed a core, but if $d$ is even then all of the connection set elements have odd weight and thus it is bipartite, a contradiction to $X$ being a core.\qedn

We can now prove that the only graphs meeting the degree bound with equality are the folded cubes of odd order:

\begin{theorem}\label{thm:folded}
Suppose that $X$ is a core of a cubelike graph, has degree $d \ge 2$, and has $2^{d-1}$ vertices. Then $d$ is odd and $X$ is the folded cube of order $d$.
\end{theorem}
\proof
%Suppose that $X^\bullet$, $Y$, $Z$, $a_1, \ldots a_k$, and $r$ are as in Theorem~\ref{thm:cover}, and that $X^\bullet$ further satisfies the hypotheses of this theorem.

Let $e_1, \ldots, e_d$ be the standard basis vectors of $\mathbb{Z}^d$, and let $\varphi$ be the covering map from the $d$-cube to $X$ guaranteed by Theorem~\ref{thm:cover}. It will suffice to show that $X$ is cubelike and then we can apply Lemma~\ref{lem:cubecores}. To do this we will show that the fibres of $\varphi$ are all cosets of a two element subgroup.

%First, if $a_1, \ldots, a_k$ are not linearly independent, then $Y$ has at most $2^{k-1}$ vertices and thus $r\restriction_Y$ is an isomorphism and therefore $X^\bullet$ is cubelike. So we may assume that $a_1, \ldots, a_k$ are linearly independent and therefore the fibres of $r\restriction_Y$ have size two. Note that this implies that $Y$ is bipartite.

Let $x$ be a vertex of $X$. Then there exists a $y \in \mathbb{Z}^d$ such that $x = \varphi(y)$. Since $X$ is not bipartite, it contains an odd cycle, and since it must be vertex transitive the vertex $\varphi(y)$ must be contained in some shortest odd cycle. Since $\varphi$ is a covering map, the vertices of this shortest odd cycle can be written as
\[\varphi(y), \varphi(y+e_{i(1)}), \varphi(y+e_{i(1)}+e_{i(2)}), \ldots, \varphi(y+e_{i(1)}+ \ldots + e_{i(g)}) = \varphi(y),\]
for some odd $g$ and $i(j) \in [d]$ for all $j \in [g]$. However, if $i(j) = i(j')$ for any $j \ne j'$, then removing $e_{i(j)}$ and $e_{i(j')}$ in all of the above terms in which they appear, and then removing the two redundant terms, would result in a shorter odd cycle of $X$, a contradiction. Therefore the $e_{i(j)}$ are distinct and we may assume without loss of generality that the vertices of the odd cycle have the form
\[\varphi(y), \varphi(y+e_{1}), \varphi(y+e_{1}+e_{2}), \ldots, \varphi(y+e_{1}+ \ldots + e_{g}) = \varphi(y).\]

Let $\Delta = \sum_{i=1}^g e_i \ne 0$. We have that $y + \Delta \ne y$, but $\varphi(y) = \varphi(y+\Delta)$. So $y$ and $y+\Delta$ are the two vertices in the preimage of $x$. We aim to show that the sum of two vertices in a fibre of $\varphi$ is always equal to $\Delta$, i.e.~every fibre is a coset of the subgroup $\{0,\Delta\}$. We will show this for the neighbors of $\varphi(y)$ and thus it will hold by connectivity of $X$.

Since $x$ was arbitrary, the above implies that the sum of two vertices in a fibre of $\varphi$ must have weight equal to the length of the shortest odd cycle of $X$, which is $g$. The neighbors of $\varphi(y)$ in $X$ can be written in the form $\varphi(y+e_i)$ for $i \in [d]$. Moreover, there is a bijection $f: [d] \to [d]$ such that $\varphi(y+e_i) = \varphi(y+\Delta + e_{f(i)})$ for all $i \in [d]$. We want to show that $f(i) = i$. The vertices in the preimage of $\varphi(y+e_i)$ have sum $\Delta+e_i+e_{f(i)}$. If $i \in [g]$ and $i \ne f(i) \in [g]$, then $\Delta+e_i+e_{f(i)}$ has weight $g-2$, a contradiction. Similarly, if $i \in [d] \setminus [g]$ and $i \ne f(i) \in [d]\setminus [g]$, then $\Delta+e_i+e_{f(i)}$ has weight $g+2$, a contradiction. So $f$ either fixes an element $i \in [d]$ or maps it to whichever of $[g]$ and $[d]\setminus [g]$ that does not contain $i$. Therefore, showing that $f$ fixes the elements of $[g]$ implies that it fixes all of the elements of $[d]$.

Now consider moving around the same shortest odd cycle we considered above, but this time we write the vertices as
\[\varphi(y+\Delta), \varphi(y+\Delta + b_1), \varphi(y+\Delta+b_{1}+b_{2}), \ldots, \varphi(y+\Delta+b_{1}+ \ldots + b_{g}) = \varphi(y+\Delta),\]
where the $b_j$ are some distinct elements of $\{e_1, \ldots, e_d\}$. As above, we have that $y+\Delta+b_{1}+ \ldots + b_{g} \ne y+\Delta$, but they have the same image under $\varphi$. This implies that
\[y = y+\Delta+b_{1}+ \ldots + b_{g} = y + \sum_{i=1}^g e_i + \sum_{i=1}^g b_i.\]
From this we see that
\[\sum_{i=1}^g e_i = \sum_{i=1}^g b_i.\]
This implies that $\{b_1, \ldots b_g\} = \{e_1, \ldots e_g\}$. Since $\varphi(y+e_1) = \varphi(y+\Delta + b_1)$ we have that $e_{f(1)} = b_1$, and since $b_1 \in \{e_1, \ldots, e_g\}$, we have that $f(1) \in [g]$. By the above this implies that $f(1) = 1$. By permuting the elements of $\{e_1, \ldots, e_g\}$, we can obtain different shortest odd cycles containing $\varphi(y)$ which will similarly imply that $f(i) = i$ for all $i \in [g]$. From this and the above arguments we see that $f$ fixes every element of $[d]$ and therefore the sum of two vertices in the preimage of any neighbor of $\varphi(y)$ is equal to $\Delta$, and by connectivity this is true for every fibre of $\varphi$.

By the above and Lemma~\ref{lem:cosets}, we have that $X$ is cubelike, and then by Lemma~\ref{lem:cubecores} we have proven the theorem.\qedn

Together, Corollary~\ref{cor:degreebound} and Theorem~\ref{thm:folded} immediately give the following corollary:

\begin{corollary}\label{cor:even}
If $X$ is the core of a cubelike graph and has even degree $d$, then $|V(X)| \le 2^{d-2}$.
\end{corollary}

Letting $d = 4$ in the above corollary, we see that any degree 4 core of a cubelike graph has at most 4 vertices, but no such (simple) graph exists. It is also not too difficult to use Theorem~\ref{thm:folded} to show that the only degree 5 core of a cubelike graph is the Clebsch graph.

\section{Small Cores}\label{sec:fewverts}

In this section we combine all of our tools from previous sections to show that if $X$ is the core of a cubelike graph and $|V(X)| \leqslant 32$, then $X$ is itself cubelike.  In the next subsection, we prove this by hand when $|V(X)| \leqslant 16$, and in the following we use a computer to prove the same for $|V(X)|=32$.

\subsection{A fistful of vertices $\ldots$}

\begin{theorem}\label{thm:fewverts}
Suppose that $X$ is the core of a cubelike graph and $|V(X)| \le 16$. Then $X$ is either $K_2, K_4, K_8, K_{16}$, or the Clebsch graph or its complement.
\end{theorem}
\proof
We will start with the $|V(X)| = 16$ case. We know that the core of an abelian Cayley graph is complete if and only if the clique-coclique bound holds with equality. Furthermore, the clique-coclique bound holds with equality for an abelian Cayley graph if and only if it holds for its core. Therefore, the clique-coclique bound does not hold for the core of a cubelike graph unless that core is a complete graph.

Recall that a cubelike graph cannot have clique number equal to three, and therefore neither can the core of a cubelike graph. Thus once we prove the core of a cubelike graph contains a $K_3$, it must contain a $K_4$. This also holds for the complement of the core of a cubelike graph by Corollary~\ref{cor:cayorbits}.

Suppose that $X$ is the core of a cubelike graph with $|V(X)| = 16$. If $X$ and its complement contain a $K_3$, then they both contain a $K_4$ and thus the clique-coclique bound holds with equality, a contradiction. So either $X$ or its complement is triangle-free.

Suppose that $X$ is triangle-free. Since it is a core, it must not be bipartite. It therefore has odd girth $2g+1 \ge 5$ and thus must contain $\mathbb{Z}_2[C_{2g+1}]$ as an induced subgraph by Lemma~\ref{lem:FQsub}. If $2g+1 \ge 7$ this graph has more than 16 vertices, a contradiction. Therefore, $2g+1 = 5$ and $X$ contains an induced Clebsch graph which has 16 vertices and so $X$ must be the Clebsch graph.

Now suppose that the complement of $X$ is triangle-free. If the complement is empty, then $X = K_{16}$. If the complement is a non-empty bipartite graph, then it must have a perfect matching since it is vertex transitive, and therefore we can color $X$ with 8 colors. Furthermore, one side of the bipartition of the complement of $X$ will induce a $K_8$ in $X$. This implies that $X$ is homomorphically equivalent to $K_8$, a contradiction. Otherwise the complement has an odd cycle and since it is homomorphically equivalent to a cubelike graph by Corollary~\ref{cor:cayorbits}, it must contain a folded cube as an induced subgraph. By a similar argument to the previous case this folded cube must be exactly the Clebsch graph. Therefore $X$ is the complement of the Clebsch graph. These are the only possibilities for $|V(X)| = 16$.

Suppose $|V(X)| = 8$. If $X$ contains a triangle then it contains a $K_4$. If its complement is not empty the clique-coclique bound will hold with equality which is a contradiction. Therefore if $X$ has a triangle it must be $K_8$. If $X$ is triangle-free it must contain an induced folded cube with at least 16 vertices, a contradiction. Therefore $K_8$ is the only core of a cubelike graph on 8 vertices.

It is easy to see that the only cores on four and two vertices are $K_4$ and $K_2$ respectively.\qedn

\subsection{$\ldots$ and a few vertices more}\label{subsec:32vertices}

For the case where the putative core of a cubelike graph has 32 vertices, we require the aid of a computer, and need to use all of our previous results to rule out every non-cubelike vertex-transitive graph as a possible core.

In this section, we describe computations showing that every vertex-transitive graph on $32$ vertices that meets our necessary conditions to be the
core of a cubelike graph is actually cubelike itself. Thus if $Z$ is a counterexample to Conjecture~\ref{conj:cubecore}, namely a cubelike graph $Z$ with a non-cubelike core $Z^\bullet$, then $|V(Z^\bullet)| \geq 64$.

We start with the complete list of the transitive permutation groups of degree $32$, which was found by Cannon \& Holt \cite{cannonholt} and is available in the computer algebra system {\sc Magma}. From this, it is relatively straightforward to compute all of the  $677116$ connected vertex-transitive graphs of order $32$, and this list forms our initial set of candidates for the core of a (possibly much larger) cubelike graph.  From this list, we remove graphs that cannot be the core of any cubelike graph, because they violate one or more of the conditions outlined in the lemmas of the previous sections. In effect each of these lemmas is used as a \emph{filter} to eliminate unsuitable graphs from the list of candidates. After applying what is (in retrospect) a surprising number and variety of these filters, every possibility is eliminated.

%Of course this approach to proving the conjecture will only work if enough properties of the core of a cubelike graph can be found that only an actual cubelike graph can satisfy them all. 

There is a large number and variety of vertex-transitive graphs on $32$ vertices, and so we feel that our results provide considerable supporting evidence in favour of the conjecture. Proving the analogous result for vertex-transitive graphs on $64$ vertices, which would provide even stronger evidence, is not possible using these techniques. The sheer number of groups and graphs will certainly not be manageable without considerably stronger theoretical results constraining the structure of the core of a cubelike graph. Ultimately, to prove the conjecture in its entirety along these lines, it will be necessary to find such strong structural properties of the core of a cubelike graph that they are sufficient to characterise cubelike graphs.

\begin{enumerate}

\item Integer Eigenvalues: from $677116$ graphs to $8648$ graphs

By Theorem~\ref{thm:cover}, if the core of a cubelike graph is $d$-regular, then it is covered by the $d$-cube $Q_d$, and hence its eigenvalues are a submultiset of the eigenvalues of $Q_d$ as noted in Corollary~\ref{cor:spectra}. As $Q_d$ has all integer eigenvalues, so does the core of any cubelike graph. We used {\sc Magma} to compute and factor the characteristic polynomial of each of the candidate graphs, retaining only those with integral spectrum.

This proved to be a very powerful test, eliminating just over $98.72\%$ of the candidate graphs, leaving only $8648$ to proceed to the next testing phase.

\item Clique-Coclique Bound Equality: from $8648$ graphs to $1966$ graphs

By Lemma~\ref{lem:cliquecoclique}, if $X$ is vertex transitive and homomorphically equivalent to a cubelike graph, then it satisfies the clique-coclique bound with equality if and only if $X^\bullet$ is complete. Therefore, if $X$ is the core of a cubelike graph that satisfies the clique-coclique bound with equality, then it must itself be a complete graph. In this case, $X$ is cubelike since it must have a power of two vertices. Thus no graph $X$ with $\alpha(X)\omega(X) = |V(X)|$ can be a non-cubelike core of a cubelike graph, and so these may be filtered out from the list of candidate graphs.

%If $X$ is a vertex-transitive graph on $n$ vertices, then $\alpha(X) \omega(X) \le n$. If equality holds in this bound, and $X$ is a normal Cayley graph, then $\chi(X) = \omega(X)$, and so the core of $X$ is complete.  
%So suppose then that $X$ is a graph such that $\alpha(X) \omega(X) = n(X)$ and also that is the core of a cubelike graph $Z$. As $Z$ has the same clique number and fractional chromatic number (i.e., $n/\alpha$) as $X$, it follows that $\chi(Z) = \omega(Z)$ and as $Z$ is a normal Cayley graph, it follows that the core of $Z$ is complete, and so $X$ is complete. Therefore no graph with $\alpha(X) \omega(X) = n$ can be non-cubelike and also the core of a cubelike graph, and so these may be filtered out from the list of candidate graphs.

\item Generously Transitive: from $1966$ to $318$ graphs

As a cubelike graph is generously transitive, its core must also be generously transitive and so we can remove any graph from the list that does not have a generously transitive automorphism group.  After this test, the list now contains only $318$ candidate cores.

\item Not Cubelike: from $318$ to $196$ graphs

There are only $1372$ cubelike graphs on $32$ vertices, and so it is now easy to decide which of the remaining candidates are cubelike by 
checking whether they are isomorphic to any graph in the list
of $1372$ graphs. 
This test reduced the number of candidate graphs from $318$ to $196$.

\item Cliques of Orbital Graphs: from $196$ to $32$ graphs

If $X$ is the core of a cubelike graph, then by Corollary~\ref{cor:cayorbits} each of its orbital graphs is homomorphically equivalent to a cubelike graph. As homomorphically equivalent graphs have equal clique number, it follows that no orbital graph of $X$ can have clique number exactly three. Although $X$ may have many orbital graphs, it is not usually necessary to construct all of them. In particular, if \emph{any} orbital graph of $X$ has clique number $3$, then there is an orbital graph of $X$ whose edge set is the union of at most three orbitals with the same property (just take the union of the orbitals containing each of the three edges of the triangle).  Most of the $196$ candidates have an orbital graph with clique number $3$, and this brings our list down to $32$ graphs.

\item Eigenvalues of Orbital Graphs: from $32$ to $20$ graphs

If $X$ is the core of a cubelike graph, then by Corollary~\ref{cor:spectra} the spectrum of a connected component of any $d$-regular orbital graph of $X$ is a submultiset of the spectrum of the cube $Q_d$. In this test, not only do we check that the eigenvalues are integers, but also that the multiplicity of each eigenvalue is no greater than the multiplicity of that eigenvalue in the spectrum of $Q_d$.

%[Thorem 4.2] any $d$-regular orbital graph of $X$ is covered by the $d$-cube and so its spectrum is a submultiset of the spectrum of the cube $Q_d$. In this test, not only do we check that the eigenvalues are integers, but also that the multiplicity of each eigenvalue is no greater than the multiplicity of the that eigenvalue in the spectrum of $Q_d$.

\item Not a Core: from $20$ to $18$ graphs

Testing whether a graph is a core is a difficult computational task, mostly because showing that a graph actually \emph{is} a core involves an exhaustive search to demonstrate that the graph has no non-surjective endomorphisms. However as there are so few remaining graphs and they are of modest size, the {\tt digraphs} package in {\sc Gap} can easily determine that $18$ of the $20$ are cores, while eliminating $2$ non-cores that have passed all previous tests.

\item Hom-Idempotence: from $18$ to $4$ graphs

%By Corollary~\ref{cor:cayorbits}, any orbital graph of the core of a cubelike graph is homomorphically equivalent to a cubelike graph. Since cubelike graphs are normal Cayley graphs, combining this with Lemma~\ref{lem:homidem} tells us that any orbital graph $X$ of the core of a cubelike graph must be hom-idempotent and therefore homomorphically equivalent to $\text{Sh}(X^\bullet)$. Since $\text{Sh}(X^\bullet) \to X^\bullet \to X$ always holds, the latter is equivalent to $X \to $\text{Sh}(X^\bullet)$. Therefore, if a 

Since any cubelike graph is a normal Cayley graph, by Lemma~\ref{lem:homidem} we have that the core $X$ of a cubelike graph must be hom-idempotent, and moreover must be homomorphically equivalent to (and therefore the core of) its shift graph $\text{Sh}(X)$. As $\text{Sh}(X) \to X$ always holds, the latter is equivalent to the existence of a homomorphism from $X$ to $\text{Sh}(X)$.

The $18$ remaining graphs have automorphism groups of orders ranging from $256$ to $1536$ and so it is not difficult to compute the shift graphs. For many of these shift graphs, it is then possible to use the  {\tt digraphs} package to determine whether or not the graphs are hom-idempotent. For some of the graphs though, this computation takes too long. However, employing Corollary~\ref{cor:cayorbits} we see that any orbital graph of the core of a cubelike graph must be hom-idempotent. Thus if a candidate graph has an orbital graph that is not hom-idempotent, it can be removed from the list of candidates.

After this stage, only $4$ graphs remain; these have not been eliminated because they \emph{are} hom-idempotent, as they are all Cayley graphs for the abelian group $Z_4 \times Z_8$.

\item Cubelike Hulls: from $4$ to zero graphs

The graph $X = {\rm Cay}(Z_4 \times Z_8, \{  (1,0), (0,6), (0,3), (0,7), (1,5), (1,1), (1,6), (2,2) \})$ from our list has $\omega(X) = 5$. By Lemma~\ref{lem:Z2ext}, if $X$ is homomorphically equivalent to a cubelike graph then we must that $\ZZ{K_5} \to X$. However it is easy to check (using {\tt digraphs}) that there is no homomorphism from $\ZZ{K_5}$ to $X$, and therefore $X$ is not homomorphically equivalent to a cubelike graph. The remaining three graphs all have $X$ as an orbital graph and hence none of them are the core of a cubelike graph by Corollary~\ref{cor:cayorbits}.

%The graph $X = {\rm Cay}(Z_4 \times Z_8, \{  (1,0), (0,6), (0,3), (0,7), (1,5), (1,1), (1,6), (2,2) \})$ has $\omega(X) = 5$. If $X$ were homomorphically equivalent to a cubelike graph $Z$, then $\omega(Z) = \omega(X) = 5$ and so $\mathbb{Z}_2(5) \rightarrow Z \rightarrow X$. However it is easy to check (using {\tt digraphs}) that there is no homomorphism from $Z_2(5)$ to $X$, and therefore $X$ is not homomorphically equivalent to a cubelike graph. The remaining three graphs all have $X$ as an orbital graph and hence none of them are homomorphically equivalent to a cubelike graph either. 

\end{enumerate}

At this point, we have shown the following:

\begin{theorem}\label{thm:32}
If $X$ is the core of a cubelike graph and $|V(X)| \le 32$, then $X$ is cubelike.
\end{theorem}

Applying the Degree Bound, we also obtain the following:

\begin{corollary}
If $X$ is the core of a cubelike graph and $X$ has valency at most 7, then $X$ is cubelike.
\end{corollary}

\section{Discussion}

In this work we have shown that both the core of a cubelike graph, and the orbital graphs of this core, must share a variety of graph-theoretical properties. In particular, we have the following:

\begin{theorem}
Suppose that $X$ is the core of a cubelike graph, and let $Y$ be a connected component of an orbital graph of $X$ with valency $d$. Then the following hold:
\begin{enumerate}
\item $Y$ is homomorphically equivalent to a cubelike graph. In particular, $\omega(Y),\chi(Y) \ne 3$.
\item $Y$ is generously transitive.
\item There is a covering map from the $d$-cube $Q_d$ to $Y$. Thus $Y$ has at most $2^d$ vertices.
\item If $Y$ is a core other than $K_2$, then $|V(Y)| \le 2^{d-1}$ with equality if and only if $d$ is odd and $Y$ is the folded cube of order $d$.
\item The eigenvalues of $Y$ are a sub-multiset of the eigenvalues of $Q_d$ and are thus integers.
\item If $\alpha(Y)\omega(Y) = |V(Y)|$ then the core of $Y$ must be complete.
\item $Y$ is hom-idempotent. If $Y$ is a core then this is equivalent to $Y$ being a subgraph of its shift graph.
\item If $Z \to Y$, then $\ZZ{Z} \to Y$ for any graph $Z$.
\item If $Y$ has odd girth $g$, then $Y$ contains the folded cube of order $g$ as an induced subgraph.
\end{enumerate}
\end{theorem}

As we have seen, the combination of these properties is very restrictive, allowing us to rule out all non-cubelike vertex-transitive graphs on up to $32$ vertices as possible cores of cubelike graphs. While we are unable to show that this list of properties suffices to characterize cubelike cores, we also know of no non-cubelike core that satisfies all of them, and so we cannot definitely say that the list does \emph{not} characterise cubelike cores. The next place to look for such a graph is amongst the vertex-transitive graphs on $64$ vertices. However, the vertex-transitive graphs on $64$ vertices are not known, and even if they were, it is likely that there would be vast numbers of them, probably enough to overwhelm the required tests. It seems that a new approach is needed, even to rule out the $64$-vertex cores. In the next section we outline a possible strategy which involves working directly with transitive groups, on the grounds that each group will have numerous orbital graphs that can, in some circumstances, be dealt with \emph{en masse}.

\subsection{A possible approach for 64 vertices}

One approach to the construction of vertex-transitive groups is to perform a breadth-first traversal of the subgroup lattice of a fixed group (usually a particular wreath product). At each stage, a group of maximum order is removed from the queue, its maximal transitive subgroups are computed, and those not conjugate to groups already in the queue are added to the queue. The orbital graphs of a  group $\Gamma$ form a subset of the orbital graphs of any subgroup of $\Gamma$, and so in certain circumstances an entire branch of the subgroup lattice can be pruned.  

More precisely, at each step we process a transitive permutation group $\Gamma$, first constructing all of its orbital graphs. Depending on $\Gamma$, exactly one of the following will hold:
\begin{enumerate}
\item Every orbital graph of $\Gamma$ is cubelike.
\item Some orbital graphs of $\Gamma$ are not cubelike but they have cubelike cores.
\item All of the orbital graphs of $\Gamma$ that are cores are cubelike, but there are some that have (strictly smaller) non-cubelike cores.
\item There is some orbital graph of $\Gamma$ that is both a core and non-cubelike.
\end{enumerate}
Depending on which case we are in, we will either continue branching, terminate, or stop with an inconclusive result, storing the group $\Gamma$ for later analysis. In the first two cases we would continue branching\footnote{In the second case, it is possible that we could terminate the branching. This could happen if $\Gamma$ has a non-cubelike orbital graph $X$ such that the core of $X$ is cubelike, but $X$ does not have integer eigenvalues for instance.}, while in the third case, we would terminate the branching. In the final case, we would either terminate if we found an orbital graph which failed one of our tests, or we would have a potential counterexample which we would save and keep for further analysis, stopping the branching.

%To understand why we are able to terminate in the cases specified above, suppose that $Y$ is a counterexample to the conjecture: a non-cubelike graph that is the core of a cubelike graph. We claim that the graph $Y$ (or some other counterexmple) will be found by our algorithm. First, note that 

To understand why we are able to terminate in the cases specified above, suppose that the permutation group $\Gamma$ currently being considered has an orbital graph $X$ that fails one of our tests. Then we know that any other graph that has $X$ as an orbital graph cannot be a counterexample to the conjecture, and so we may immediately rule it out. What are the graphs that have $X$ as an orbital graph? Those whose automorphism group is a subgroup of $\aut(X)$. Since $X$ was an orbital graph of $\Gamma$, we have that $\Gamma \le \aut(X)$. So if $Y$ is a counterexample to the conjecture, then we cannot have that $\aut(Y) \le \Gamma$. Thus $\aut(Y)$ will be reached by some other branch of the algorithm (or that branch of the algorithm will stop with some other potential counterexample before reaching $\aut(Y)$).

%$\Gamma' \le \Gamma$ and $Y$ is an orbital graph of $\Gamma'$ that is a countere, then either $\aut(Y) \le \Gamma$ no orbital graph of any orbital graph contains an isomorphic copy of $\mathbb{Z}_2^6$ as a subgroup, then every orbital graph of $\Gamma$ is cubelike. 

In the third case, we can terminate immediately because if one of the orbital graphs $X$ has a strictly smaller core that is not cubelike, then this core cannot be the core of any cubelike graph by our results for 32 or fewer vertices. Thus $X$ is not homomorphically equivalent to a cubelike graph, and so cannot be the orbital graph of any counterexample. In the final case, we may be unlucky and find a graph $Y$ that is a non-cubelike core but we cannot decide whether it is the core of some cubelike graph or not. This would happen if it passed all of the tests we could reasonably perform on it. It may be that $Y$ is not the core of any cubelike graph, but it is hom-idempotent, is covered by a $d$-cube, is generously transitive, does not have chromatic number three, and $Z \to Y \Rightarrow \ZZ{Z} \to Y$ for all graphs $Z$ up to say 20 vertices (and this holds for all non-cubelike orbital graphs of $Y$ as well). The fundamental problem seems to be that even if $Y$ were a legitimate counterexample, we do not have any efficient way of showing this. Either we prove that $\ZZ{Y} \to Y$, which is computationally impossible for $|V(Y)| = 64$, or we somehow construct some cubelike graph and prove through other means that its core is $Y$. On the other hand, it is possible that we do not run into this issue at all, which was the case for 32 vertices.

Note that if the group $\Gamma$ contains a regular elementary abelian 2-group, then we will necessarily be in the first case. So our first step is to find all of the maximal permutation groups of degree 64 without such a subgroup. However, even this seems to be beyond our current computational limits.

The advantage that this approach does have is that when we terminate one of the branches, we are essentially excluding many graphs/permutation groups that we never had to consider directly, and so the total number of graphs/groups we have to consider may be significantly more manageable than if we simply construct all of the transitive graphs on 64 vertices. However, another difficulty is that it is possible for a permutation group to be visited by many different branches, and this may cost us more than we have gained. Finding a way to avoid or minimize these revisits would be a crucial step towards making this approach viable. But for now, tackling the 64 vertex case seems to be out of reach.

\subsection{An idea towards proving the conjecture}

Though we do not know how to prove the conjecture, we present an idea that may or may not be useful to other researchers who wish to work on this problem. This idea is related to the second false lead from Section~\ref{subsec:falseleads}. There we showed that it is not always the case that a cubelike graph $Z$ has a homomorphism to its core whose fibres are cosets of a single subgroup. Here, we show that the conjecture holds if and only if for any graph $X$ which is the core of some cubelike graph, there is a homomorphism from $\ZZ{X}$ to $X$ whose fibres are cosets of a single subgroup.

Suppose that $X$ is the core of a cubelike graph. We know from Lemma~\ref{lem:Z2core} that $X$ is the core of $\ZZ{X}$. If $X$ is cubelike, then by Lemma~\ref{lem:Z2} we have that there is a (vertex- and edge-surjective) graph homomorphism $\varphi$ from $\ZZ{X}$ to $X$ which is also a group homomorphism. Since it is a group homomorphism, the fibres of $\varphi$ are all cosets of a fixed subgroup: the kernel of $\varphi$. Conversely, if $X$ is not cubelike, then by Lemma~\ref{lem:cosets} there exists no such homomorphism from $\ZZ{X}$ to $X$. Therefore, if $X$ is a core, then $X$ is cubelike if and only if there exists a homomorphism from $\ZZ{X}$ to $X$ whose fibres are all cosets of a fixed subgroup. So to prove the conjecture, it suffices to show that for any core $X$ such that $\ZZ{X} \to X$, there exists a homomorphism from $\ZZ{X}$ to $X$ whose fibres are cosets of a fixed subgroup.

We remark that, for a core $X$ with $2^n$ vertices, any homomorphism from $\ZZ{X}$ to $X$ is also a $2^n$-coloring of the halved $2^n$-cube $\frac{1}{2}Q_{2^n}$. To see this recall that for any $u \in V(X)$ the set $\{e_u+e_v: v \in V(X)\}$ induces a copy of $X$ in $\ZZ{X}$. This implies that for any two vertices $u,v \in V(X)$, the pair $(0,e_u+e_v)$ is contained in a copy of $X$ in $\ZZ{X}$. In the case where there exists some homomorphism from $\ZZ{X}$ to $X$, the latter is the core of the former. Therefore, it follows from Corollary~\ref{cor:dist2} that no endomorphism of $\ZZ{X}$ can identify any two vertices of whose sum has weight two. The pairs of vertices whose sum has weight two are exactly the edges of the halved cube, and thus the fibres of any homomorphism from $\ZZ{X}$ to $X$ must be independent sets in the halved $2^n$-cube $\frac{1}{2}Q_{2^n}$. Since $X$ has $2^n$ vertices, this gives a partition of the vertices of the halved cube into $2^n$ independent sets, i.e., a $2^n$-coloring of $\frac{1}{2}Q_{2^n}$.

One is now tempted to prove that the color classes of any $2^n$-coloring of $\frac{1}{2}Q_{2^n}$ are cosets of some fixed subgroup. Unfortunately, this is not the case. Recall the $2^n$-coloring $\varphi$ of $\frac{1}{2}Q_{2^n}$ described in Section~\ref{subsec:falseleads}: $\varphi$ was the linear extension of a map $f$ from the standard basis vectors of $\mathbb{Z}_2^{2^n - 1}$ to the nonzero elements of $\mathbb{Z}_2^n$. The fibres of this map were cosets of the kernel of $\varphi$, but by choosing $f$ differently, we can obtain a different kernel. If the subgroup of $\mathbb{Z}_2^{2^n - 1}$ generated by both of these kernels is not the whole group\footnote{We leave it to the reader to verify that this is possible for large enough $n$.}, then we can partition this subgroup with cosets of the first kernel, and partition the remainder of the group with cosets of the second kernel. This gives a $2^n$ coloring whose color classes are all cosets of subgroups, but not the same subgroup. 

%We assigned the nonzero elements of $\mathbb{Z}_2^n$ to the $2^n - 1$ coordinates of the elements of $\mathbb{Z}_2^{2^n - 1}$. The color classes were the cosets of the subgroup of $\mathbb{Z}_2^{2^n - 1}$ consisting of the elements 

However, one only needs to show that there \emph{exists} some homomorphism from $\ZZ{X}$ to $X$ whose fibres are cosets of some fixed subgroup in order to show that $X$ is cubelike. So perhaps this idea is still useful.

\subsection*{Acknowledgements}

The authors wish to thank the members of the Centre for the Mathematics of Symmetry and Computation (UWA) for their interest when this problem was presented at the CMSC annual research retreat, and particularly Gabriel Verret for several subsequent discussions. Part of this work was completed during a visit of DR to the Centre for the Mathematics of Symmetry and Computation (UWA) which was funded by the Cheryl E. Praeger Visiting Research Fellowship.

%\bibliographystyle{plainurl}
%
%\bibliography{CoresofCubelikeGraphs.bbl}

\end{document}